\documentclass[12pt]{article}
\usepackage{mathrsfs}
\usepackage{stmaryrd}
\usepackage{amsfonts,amsmath,amssymb,amscd}
\usepackage[vcentermath,enableskew]{youngtab}
\usepackage{dsfont}
\usepackage{shadow}
\usepackage{graphicx}
\usepackage{color}
\usepackage{longtable}
\usepackage{pstricks,multido}
\usepackage{hyperref}
\usepackage{qtree}
\allowdisplaybreaks

\newtheorem{theorem}{Theorem}[section]
\newtheorem{lemma}[theorem]{Lemma}

\newtheorem{corollary}[theorem]{Corollary}

\newtheorem{defi}[theorem]{Definition}

\def\qed{\hfill \rule{4pt}{7pt}}

\def\pf{\noindent {\it{Proof.} \hskip 2pt}}

\def\x{\mathbf{x}}
\def\f{\mathbf{f}}
\def\y{\mathbf{y}}
\def\S{\mathfrak{S}}
\def\G{\mathfrak{G}}
\def\SRT{{\rm SRT}}
\def\SVRT{{\rm SVRT}}
\def\CSBL{{\rm CSBL}}

\numberwithin{equation}{section}
\numberwithin{figure}{section}

\begin{document}

\begin{center}
{\large {\bf Set-valued  Rothe
Tableaux and Grothendieck Polynomials}}

\vskip 4mm
{\small  Neil J.Y. Fan$^1$ and Peter L. Guo$^2$}

\vskip 4mm

$^1$Department of Mathematics\\
Sichuan University, Chengdu,
Sichuan 610064, P.R. China
\\[3mm]

$^{2}$Center for Combinatorics, LPMC\\
Nankai University,
Tianjin 300071,
P.R. China

\vskip 4mm

$^1$fan@scu.edu.cn,
$^2$lguo@nankai.edu.cn

\end{center}

\begin{abstract}
The notion of   set-valued Young  tableaux was introduced by Buch in his study of the Littlewood-Richardson rule for stable Grothendieck polynomials. Knutson, Miller and Yong
showed that the double Grothendieck polynomials of   2143-avoiding  permutations  can be generated by
 set-valued Young tableaux.
In this paper, we introduce the structure  of set-valued Rothe tableaux of permutations. Given the Rothe diagram $D(w)$ of a permutation $w$, a  set-valued Rothe tableau of shape $D(w)$ is
a filling of finite nonempty subsets of positive integers into the squares of $D(w)$ such that the rows are weakly decreasing and the columns are strictly increasing. We show that the double Grothendieck polynomials of 1432-avoiding permutations can be generated by    set-valued Rothe tableaux. When restricted to  321-avoiding permutations, our formula specializes to the tableau  formula for double Grothendieck polynomials  due to Matsumura.
Employing the properties of tableau complexes given  by Knutson, Miller and Yong, we obtain  two alternative  tableau formulas for the double Grothendieck polynomials of    1432-avoiding permutations.

\end{abstract}

\section{Introduction}

Let $S_n$ denote the set of permutations on $\{1,2,\ldots,n\}$.
The double Grothendieck polynomials $\G_{w}(\x,\y)$ indexed by permutations $w\in S_n$ were introduced by Lascoux and Sch\"{u}tzenberger \cite{LaSh} as polynomial representatives of the equivariant $K$-theory classes of structure sheaves of Schubert varieties in the flag manifold.
These polynomials were originally defined based on the isobaric
divided difference operators. Several combinatorial models have  been developed  to generate $\G_{w}(\x,\y)$, see, for example, \cite{BR,FoKi,KnMi,KnMi2,Le}.

On the other hand, tableau formulas for $\mathfrak{G}_{w}(\x,\y)$ have been found for  specific families of  permutations.
Based on the algebraic geometry of matrix Schubert varieties, Knutson, Miller and Yong \cite{KnMiYo-1} showed that for a 2143-avoiding permutation $w$ (also called a vexillary permutation), $\mathfrak{G}_{w}(\x,\y)$  can be generated by flagged  set-valued   Young tableaux.
A permutation $w=w_1w_2\cdots w_n\in S_n$ is 2143-avoiding if there do not exist indices
$1\le i_1<i_2<i_3<i_4\le n$ such that $w_{i_2}<w_{i_1}<w_{i_4}<w_{i_3}$.
 Set-valued Young tableaux  were introduced by Buch
\cite{Bu} in his study of the Littlewood-Richardson rule for stable Grothendieck polynomials.
Restricting to
semistandard Young tableaux (namely, set-valued Young tableaux with each set containing
a single integer), the Knutson-Miller-Yong  formula
specializes to the tableau formula for  the Schubert polynomial  $\mathfrak{S}_{w}(\x)$ of 2143-avoiding permutations  due to Wachs \cite{Wa}.

By introducing the structure of tableau  complexes and  utilizing the tools of commutative algebra, Knutson, Miller and Yong \cite{KnMiYo-2} found two other tableau formulas of $\G_w(\x,\y)$ for 2143-avoiding permutations
in terms of semistandard Young tableaux and limit set-valued Young tableaux,
respectively. A limit set-valued Young tableau is an assignment of finite  nonempty subsets of positive
integers to the squares of a Young diagram such that one can pick out an integer from
each square to form a semistandard Young tableau.

Recently, Matsumura \cite{Matsumura-2} provided  a tableau formula of $\G_w(\x,\y)$ for 321-avoiding permutations. A permutation $w=w_1w_2\cdots w_n$ is 321-avoiding if there do not exist indices
$i_1<i_2<i_3$ such that $w_{i_1}>w_{i_2}>w_{i_3}$.
To a 321-avoiding permutation $w$, one can associate
a skew Young diagram, denoted  $\sigma(w)$. Matsumura \cite{Matsumura-2} showed that for a 321-avoiding permutation $w$,
$\mathfrak{G}_{w}(\x,\y)$ can be generated by flagged  set-valued  tableaux of
shape $\sigma(w)$. This formula generalizes
the tableau formula for the single Grothendieck polynomial $\G_w(\x)$  of a 321-avoiding permutation
given by
Anderson,  Chen and Tarasca \cite{Anderson}.
When restricted  to semistandard Young tableaux,
it specializes to the formula for the double Schubert polynomial $\S_w(\x,\y)$
of a 321-avoiding permutation obtained by Chen, Yan and
Yang \cite{Chen}.

In this paper, we introduce the structure  of  set-valued Rothe tableaux.
Let $D(w)$ be the Rothe diagram of a permutation $w$.
A  set-valued Rothe tableau of shape $D(w)$  is a filling of finite nonempty
subsets of positive integers into the squares of $D(w)$  such that
the sets in each row are {\it weakly decreasing} and the sets in each column are
{\it strictly increasing}. As defined by Buch \cite{Bu}, for two finite nonempty sets $A$ and $B$ of positive integers, $A<B$ if $\max A<\min B$, and $A\leq B$ if $\max A\leq  \min B$.
It was noticed by Billey, Jockusch and  Stanley \cite{BiSt} that  when $w$ is a 321-avoiding permutation,   $D(w)$ is a skew Young diagram
after a reflection about a vertical line. In this case, each row in a set-valued Rothe tableau of
shape $D(w)$ is weakly increasing after a reflection about a vertical line, and thus  becomes a set-valued Young tableau.
 Hence  set-valued Rothe tableaux can be viewed as a generalization of
  set-valued Young tableaux from Young diagrams to Rothe diagrams.

Our main objective   is  to establish  set-valued Rothe tableau
formulas of $\G_w(\x,\y)$ for a new family of permutations, namely, 1432-avoiding  permutations. A permutation $w=w_1w_2\cdots w_n$ is 1432-avoiding if there do not exist indices
$i_1<i_2<i_3<i_4$ such that $w_{i_1}<w_{i_4}<w_{i_3}<w_{i_2}$.
When restricted to 321-avoiding permutations,  we show that one of our formulas coincides with
the  formula  of  Matsumura \cite{Matsumura-2}.
It should be noted   that Stankova \cite{St} proved   that the number of 1432-avoiding permutations in $S_n$
is equal to the number of 2143-avoiding permutations in $S_n$.

\begin{figure}[h]
\setlength{\unitlength}{0.5mm}
\begin{center}
\begin{picture}(260,65)


\qbezier[60](0,0)(30,0)(60,0)\qbezier[60](0,10)(30,10)(60,10)
\qbezier[60](0,20)(30,20)(60,20)\qbezier[60](0,30)(30,30)(60,30)
\qbezier[60](0,40)(30,40)(60,40)\qbezier[60](0,50)(30,50)(60,50)
\qbezier[60](0,60)(30,60)(60,60)

\qbezier[60](0,0)(0,30)(0,60)\qbezier[60](10,0)(10,30)(10,60)
\qbezier[60](20,0)(20,30)(20,60)\qbezier[60](30,0)(30,30)(30,60)
\qbezier[60](40,0)(40,30)(40,60)\qbezier[60](50,0)(50,30)(50,60)
\qbezier[60](60,0)(60,30)(60,60)

\put(55,35){\circle*{3}}\put(45,5){\circle*{3}}
\put(25,25){\circle*{3}}\put(5,15){\circle*{3}}
\put(35,55){\circle*{3}}
\put(15,45){\circle*{3}}

\put(0,30){\line(1,0){10}}\put(0,40){\line(1,0){10}}
\put(0,20){\line(1,0){10}}
\put(0,50){\line(1,0){30}}\put(0,60){\line(1,0){30}}
\put(0,20){\line(0,1){40}}\put(10,20){\line(0,1){40}}
\put(20,50){\line(0,1){10}}\put(30,50){\line(0,1){10}}

\put(20,30){\line(1,0){10}}\put(20,40){\line(1,0){10}}
\put(20,30){\line(0,1){10}}\put(30,30){\line(0,1){10}}

\put(40,30){\line(1,0){10}}\put(40,40){\line(1,0){10}}
\put(40,30){\line(0,1){10}}\put(50,30){\line(0,1){10}}


\qbezier[60](100,0)(130,0)(160,0)
\qbezier[60](100,10)(130,10)(160,10)
\qbezier[60](100,20)(130,20)(160,20)
\qbezier[60](100,30)(130,30)(160,30)
\qbezier[60](100,40)(130,40)(160,40)
\qbezier[60](100,50)(130,50)(160,50)
\qbezier[60](100,60)(130,60)(160,60)

\qbezier[60](100,0)(100,30)(100,60)
\qbezier[60](110,0)(110,30)(110,60)
\qbezier[60](120,0)(120,30)(120,60)
\qbezier[60](130,0)(130,30)(130,60)
\qbezier[60](140,0)(140,30)(140,60)
\qbezier[60](150,0)(150,30)(150,60)
\qbezier[60](160,0)(160,30)(160,60)

\put(155,35){\circle*{3}}\put(145,5){\circle*{3}}
\put(135,25){\circle*{3}}\put(105,15){\circle*{3}}
\put(145,55){\circle*{3}}
\put(115,45){\circle*{3}}

\put(100,30){\line(1,0){10}}\put(100,40){\line(1,0){10}}
\put(100,20){\line(1,0){10}}
\put(100,50){\line(1,0){30}}\put(100,60){\line(1,0){30}}
\put(100,20){\line(0,1){40}}\put(110,20){\line(0,1){40}}
\put(120,50){\line(0,1){10}}\put(130,50){\line(0,1){10}}

\put(120,30){\line(1,0){10}}\put(120,40){\line(1,0){10}}
\put(120,30){\line(0,1){10}}\put(130,30){\line(0,1){10}}

\put(140,30){\line(1,0){10}}\put(140,40){\line(1,0){10}}
\put(140,30){\line(0,1){10}}\put(150,30){\line(0,1){10}}

\put(103,52){1}\put(103,42){2}\put(103,32){3}\put(103,22){4}
\put(113,52){1}\put(123,52){1}

\put(121,32){3} \put(125,32){2}\put(141,32){2} \put(145,32){1}


\qbezier[60](200,0)(230,0)(260,0)\qbezier[60](200,10)(230,10)(260,10)
\qbezier[60](200,20)(230,20)(260,20)\qbezier[60](200,30)(230,30)(260,30)
\qbezier[60](200,40)(230,40)(260,40)\qbezier[60](200,50)(230,50)(260,50)
\qbezier[60](200,60)(230,60)(260,60)

\qbezier[60](200,0)(200,30)(200,60)\qbezier[60](210,0)(210,30)(210,60)
\qbezier[60](220,0)(220,30)(220,60)\qbezier[60](230,0)(230,30)(230,60)
\qbezier[60](240,0)(240,30)(240,60)\qbezier[60](250,0)(250,30)(250,60)
\qbezier[60](260,0)(260,30)(260,60)

\put(255,35){\circle*{3}}\put(245,5){\circle*{3}}
\put(225,25){\circle*{3}}\put(205,15){\circle*{3}}
\put(235,55){\circle*{3}}
\put(215,45){\circle*{3}}

\put(200,30){\line(1,0){10}}\put(200,40){\line(1,0){10}}
\put(200,20){\line(1,0){10}}
\put(200,50){\line(1,0){30}}\put(200,60){\line(1,0){30}}
\put(200,20){\line(0,1){40}}\put(210,20){\line(0,1){40}}
\put(220,50){\line(0,1){10}}\put(230,50){\line(0,1){10}}

\put(220,30){\line(1,0){10}}\put(220,40){\line(1,0){10}}
\put(220,30){\line(0,1){10}}\put(230,30){\line(0,1){10}}

\put(240,30){\line(1,0){10}}\put(240,40){\line(1,0){10}}
\put(240,30){\line(0,1){10}}\put(250,30){\line(0,1){10}}

\put(203,52){\bf 1}\put(203,42){\bf 2}\put(203,32){\bf 3}\put(203,22){\bf 4}
\put(213,52){\bf 1}\put(223,52){\bf 1}

\put(221,32){\bf 2} \put(225,32){1}\put(241,32){2} \put(245,32){\bf 1}
\put(25,-10){{\small (a)}}\put(125,-10){{\small (b)}}
\put(225,-10){{\small (c)}}
\end{picture}
\end{center}
\caption{(a) The Rothe diagram $D(w)$,  (b) a set-valued Rothe tableau, (c) a limit set-valued Rothe tableau for $w=426315$.}
\label{Rothe}
\end{figure}
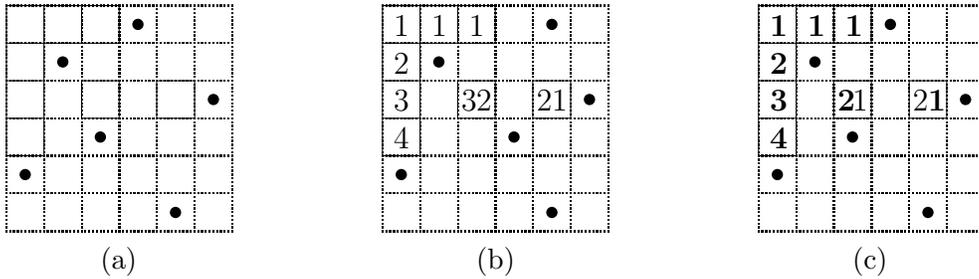

In order to state our results, we recall some definitions and notation. The Rothe diagram $D(w)$  of  a permutation $w\in S_n$ can be viewed as a geometric configuration  of the inversions of $w$.
Consider  an $n\times n$ square grid, where we use $(i,j)$ to denote the square in row $i$
and column $j$. Here the rows are numbered from top to bottom  and
the columns are numbered from left to right. For $1\leq i\leq n$,  put a dot in the square $(i, w_i)$. Then the Rothe diagram $D(w)$ consists of the squares $(i,j)$ such that there is a dot in row $i$ that is to the right of $(i,j)$, and there is a dot in column $j$ that is  below $(i,j)$.
For example, Figure \ref{Rothe}(a) is the Rothe diagram of $w=426315$.

As  aforementioned, a set-valued Rothe tableau of shape $D(w)$  is a filling of finite nonempty
subsets of positive integers into the squares of $D(w)$  such that
 the rows are weakly decreasing  and the columns are
strictly increasing. For example, Figure \ref{Rothe}(b) depicts a set-valued Rothe tableau for $w=426315$.
We say that a set-valued Rothe tableau  is flagged by a vector $\mathbf{f}=(f_1,f_2,\ldots,f_n)$ of nonnegative integers
if every  integer in row $i$  does not exceed $f_i$.
Let $\SVRT(w,\f)$ denote the set of  set-valued Rothe tableaux of shape $D(w)$ flagged by $\f$.

For a set-valued Rothe tableau $T$ and a square $B=(i,j)$ of $T$, we use $T(B)$ or $T(i,j)$ to denote the set filled in $B$.
Write  $|T|=\sum_{B\in D(w)}|T(B)|$. Let $\ell(w)$ denote the length of $w$, or equivalently, $\ell(w)=|D(w)|$. For two variables $x$ and $y$, we adopt
the following notation as used by Fomin and Kirillov \cite{FoKi}:
\[x\oplus y=x+y-xy.\]
For a square $(i,j)$ of $D(w)$, define
\begin{align}\label{statistic}
m_{ij}(w)=|\{(i,k)\in D(w)\,|\,k\leq j\}|.
\end{align}
Throughout this paper, we  use  the following specific flag \[\f_0=(1,2,\ldots,n).\]
Our main result can be stated  as follows.

\begin{theorem}\label{main}
For a permutation  $w\in S_n$,  we have the following  equivalent
statements.
\begin{itemize}
\item[(1)]  $w$ is a 1432-avoiding permutation.

\item[(2)] $\mathfrak{G}_{w}(\x,\y)$ has the following set-valued Rothe tableau formula:
\begin{align}\label{CFG-R}
\mathfrak{G}_{w}(\x,\y)=\sum_{T\in {\rm SVRT}(w,\f_0)}(-1)^{|T|-\ell(w)} \prod_{(i,j)\in D(w)}\prod_{t\in T({i,j})}(x_t\oplus y_{m_{ij}(w)+i-t}).
\end{align}
\end{itemize}
\end{theorem}

Setting $y_i=0$ in \eqref{CFG-R}, we obtain a tableau formula for
single Grothendieck polynomials  of 1432-avoiding permutations.

\begin{corollary}\label{coro}
Let $w\in S_n$ be a 1432-avoiding permutation. Then
\[
\mathfrak{G}_{w}(\x)=\sum_{T\in {\rm SVRT}(w,\f_0)}(-1)^{|T|-\ell(w)}
\prod_{(i,j)\in T}\prod_{t\in T({i,j})}x_t.
\]
\end{corollary}

The double Schubert polynomial $\mathfrak{S}_{w}(\x,\y)$
can be obtained from $\mathfrak{G}_{w}(\x,\y)$ by extracting the
monomials of  the lowest degree and then replacing  $y_i$ by $-y_i$.
Let ${\rm SRT}(w,\f)$ be the set of single-valued Rothe tableaux of shape $D(w)$
flagged by $\f$. In other words, ${\rm SRT}(w,\f)$ consists of  the Rothe tableaux   in  ${\rm SVRT}(w,\f)$ such that  the set filled in  each square  contains exactly one integer. We have the following tableau formulas for double and single Schubert
polynomials.

\begin{corollary}\label{coro}
Let $w\in S_n$ be a 1432-avoiding permutation. Then
\begin{align*}
\mathfrak{S}_{w}(\x,\y)&=\sum_{T\in {\rm SRT}(w,\f_0)}
\prod_{(i,j)\in D(w)}\prod_{t\in T({i,j})}(x_t-y_{m_{ij}(w)+i-t}),\\[5pt]
\mathfrak{S}_{w}(\x)&=\sum_{T\in {\rm SRT}(w,\f_0)}
\prod_{(i,j)\in D(w)}\prod_{t\in T({i,j})}x_t.
\end{align*}
\end{corollary}

Furthermore, by introducing the structure of Rothe tableau complexes and employing the properties of tableau complexes given  by Knutson, Miller and Yong \cite{KnMiYo-2}, we also find two alternative tableau formulas of $\G_w(\x,\y)$ for 1432-avoiding permutations. One is  given in terms of single-valued Rothe tableaux, and
the other is given in terms of limit set-valued Rothe tableaux.
A limit set-valued Rothe tableau is an assignment
of finite nonempty subsets of positive integers to the squares of a Rothe diagram such that one can pick out an integer from
each square to form a single-valued Rothe tableau.
Figure \ref{Rothe}(c) illustrates a limit set-valued Rothe tableau, where the integers in boldface form a single-valued Rothe tableau.

Let  ${\rm LSVRT}(w, \f)$   denote
the set of limit set-valued Rothe tableaux  of shape $D(w)$ flagged by $\f$.  Then we have the following two alternative tableau formulas of $\G_w(\x,\y)$ for  1432-avoiding permutations.

\begin{theorem}\label{3-1}
Let $w\in S_n$ be a 1432-avoiding permutation.
\begin{itemize}
\item[(1)]
For each square $B=(i,j)\in D(w)$, set
\[E_B=\bigcup_{T\in {\rm SRT}(w, \f_0)}T(i,j).\]
Then
\begin{align}\label{3-1-1}
\mathfrak{G}_{w}(\x,\y)= \sum_{T\in {\rm LSVRT}(w, \f_0)} \prod_{B=(i,j)\in D(w)}
\prod_{t\in T(i,j)}(x_t\oplus y_{m_{ij}(w)+i-t})\nonumber\\[5pt]
 \cdot\prod_{t\in E_B\setminus T(i,j)}(1-x_t)(1-y_{m_{ij}(w)+i-t}).
\end{align}

\item[(2)]
Given $T\in {\rm SRT}(w,\f_0)$  and a square
 $B\in D(w)$, let $Y_{T,B}$ be the set
    of positive integers $m$ such that $m$ is larger than the (unique) integer in $T(B)$ and replacing the integer in $T(B)$ by $m$ still yields a Rothe tableau in ${\rm SRT}(w, \f_0)$. Then
\begin{align}\label{3-1-2}
\mathfrak{G}_{w}(\x,\y)= \sum_{T\in {\rm SRT}(w, \f_0)}\  \prod_{B=(i,j)\in D(w)}
\prod_{t\in T(i,j)} (x_t\oplus y_{m_{ij}(w)+i-t})\nonumber\\[5pt]
    \cdot\prod_{t\in Y_{T,B}}(1-x_t)(1-y_{m_{ij}(w)+i-t}).
\end{align}
\end{itemize}
\end{theorem}

\section{Proof of Theorem \ref{main}}\label{sec3}

In this section, we aim to prove Theorem \ref{main}. For simplicity, let
\[-x \oplus y=-(x \oplus y)=-(x+y-xy).\]
Denote
 \begin{align}\label{rw}
 G_w(\x,\y)&=\sum_{T\in {\rm SVRT}(w,\f_0)}(-1)^{|T|-\ell(w)} \prod_{(i,j)\in T}\prod_{t\in T({i,j})}(x_t\oplus y_{m_{ij}(w)+i-t})\nonumber\\[6pt]
 &=(-1)^{\ell(w)}\sum_{T\in {\rm SVRT}(w,\f_0)} \prod_{(i,j)\in T}\prod_{t\in T({i,j})}(-x_t\oplus y_{m_{ij}(w)+i-t})
 \end{align}
to be the right-hand side of \eqref{CFG-R}.
We finish the proof of Theorem \ref{main} by separately proving   the following two theorems.

\begin{theorem}\label{prop1}
If $w$ is a 1432-avoiding permutation, then $\G_w(\x,\y)=G_w(\x,\y)$.
\end{theorem}

\begin{theorem}\label{prop2}
If $w$ contains a 1432 pattern, then $\G_w(\x,\y)\neq G_w(\x,\y)$.
\end{theorem}

We use the opportunity here to explain  that when  $w$ is a 321-avoiding permutation, Theorem \ref{prop1}
specializes  to the tableau formula  for $\mathfrak{G}_{w}(\x,\y)$ due to Matsumura \cite{Matsumura-2}.
To describe the tableau formula in  \cite{Matsumura-2}, let $f(w)=(f_1,f_2,\ldots,f_k)$
(respectively, $f^c(w)=(f_1^c,f_2^c,\ldots,f_{n-k}^c)$) be the increasing arrangement of the positions $i$ such that $w_i>i$ (respectively, $w_i\leq i$).
Moreover, let $h(w)=(w_{f_1},w_{f_2},\ldots,w_{f_k})$
and  $h^c(w)=(w_{f_1^c},w_{f_2^c},\ldots,w_{f_{n-k}^c})$.
It can be shown that  $w$ is 321-avoiding if and only if the sequences $h(w)$ and $h^c(w)$
are both increasing \cite{ErLi}. One may associate a skew shape $\sigma(w)=\lambda/ \mu$
to $w$ by letting
\begin{align}\label{lmu}
\lambda_i=w_{f_k}-k-(f_i-i), \ \ \ \ \ \ \mu_i=w_{f_k}-k-(w_{f_i}-i),
\end{align}
where $1\leq i\leq k$.
For a square $\alpha$ of $\sigma(w)$, let $r(\alpha)$ and $c(\alpha)$
denote the row index and the column index of $\alpha$, respectively.
\begin{corollary}[\mdseries{Matsumura \cite[Theorem 3.1]{Matsumura-2}}]
Let $w\in S_n$ be a
321-avoiding permutation. Then
\begin{align}\label{MJAP}
\mathfrak{G}_{w}(\x,\y)=\sum_{T}(-1)^{|T|-\ell(w)}
\prod_{\alpha\in \sigma(w)}\prod_{t\in T(\alpha)}(x_t\oplus y_{\lambda_{r(\alpha)}
+f_{r(\alpha)}-c(\alpha)-t+1}),
\end{align}
where $T$ ranges  over set-valued Young tableaux of shape $\sigma(w)$
flagged by $f(w)$.
\end{corollary}

\begin{figure}[h]
\setlength{\unitlength}{0.5mm}
\begin{center}
\begin{picture}(180,65)
\qbezier[60](0,0)(30,0)(60,0)\qbezier[60](0,10)(30,10)(60,10)
\qbezier[60](0,20)(30,20)(60,20)\qbezier[60](0,30)(30,30)(60,30)
\qbezier[60](0,40)(30,40)(60,40)\qbezier[60](0,50)(30,50)(60,50)
\qbezier[60](0,60)(30,60)(60,60)

\qbezier[60](0,0)(0,30)(0,60)\qbezier[60](10,0)(10,30)(10,60)
\qbezier[60](20,0)(20,30)(20,60)\qbezier[60](30,0)(30,30)(30,60)
\qbezier[60](40,0)(40,30)(40,60)\qbezier[60](50,0)(50,30)(50,60)
\qbezier[60](60,0)(60,30)(60,60)

\put(25,55){\circle*{3}}\put(5,45){\circle*{3}}
\put(15,35){\circle*{3}}\put(35,25){\circle*{3}}
\put(55,15){\circle*{3}}
\put(45,5){\circle*{3}}

\put(0,50){\line(1,0){20}}\put(0,60){\line(1,0){20}}
\put(0,50){\line(0,1){10}}\put(10,50){\line(0,1){10}}\put(20,50){\line(0,1){10}}

\put(40,10){\line(0,1){10}}\put(50,10){\line(0,1){10}}
\put(40,10){\line(1,0){10}}\put(40,20){\line(1,0){10}}

\put(-5,45){\line(1,0){70}}\put(-5,35){\line(1,0){70}}
\put(-5,25){\line(1,0){70}}\put(-5,5){\line(1,0){70}}

\put(25,-5){\line(0,1){70}}\put(55,-5){\line(0,1){70}}


\put(130,20){\line(1,0){10}}\put(130,30){\line(1,0){10}}
\put(130,20){\line(0,1){10}}\put(140,20){\line(0,1){10}}

\put(150,30){\line(1,0){20}}\put(150,40){\line(1,0){20}}
\put(150,30){\line(0,1){10}}\put(160,30){\line(0,1){10}}
\put(170,30){\line(0,1){10}}

\put(145,2){$\sigma(w)$}

\end{picture}
\end{center}
\caption{$D(w)$  and the corresponding
 skew shape $\sigma(w)$ for $w=312465$.}
\label{skew-shape}
\end{figure}
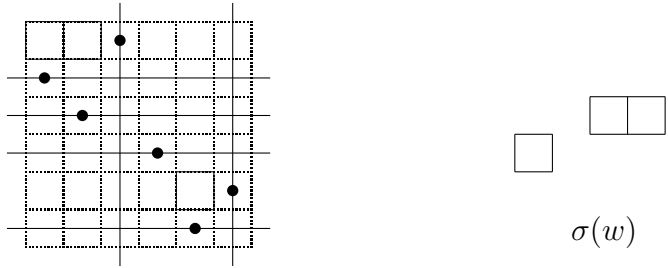

\pf
We show  that for a 321-avoiding permutation $w$, the right-hand
 side of \eqref{MJAP} is equal to  $G_w(\x, \y)$ as defined in \eqref{rw}.
As observed in \cite{BiSt}, after deleting the empty rows indexed by $f^c(w)$ and the empty columns indexed by $h(w)$ and then reflecting the resulting diagram   about  a vertical line,
 $D(w)$ coincides with the above defined
 skew shape $\sigma(w)$.
 For example, for  $w=312465$, we see that
  $f(w)=(1,5),  f^c(w)=(2,3,4,6)$ and $h(w)=(w_{f_1},w_{f_2})=(3,6)$. So the corresponding
  shew diagram $\sigma(w)$ is as illustrated in  Figure \ref{skew-shape}.

Therefore, each set-valued Rothe  tableau $T\in\SVRT(w,\f_0)$
can be viewed as a set-valued (skew) Young tableau of shape $\sigma(w)$ flagged by $f(w)$.
For a square $(i,j)\in D(w)$, assume that $\alpha$
 is the corresponding square of $\sigma(w)$. Then we need to show that
\begin{align}\label{rbr}
 \lambda_{r(\alpha)}+f_{r(\alpha)}-c(\alpha)+1=m_{ij}(w)+i.
 \end{align}
It is not hard to check that
\begin{align*}
r(\alpha)&=i-|\{t\,|\,w_t\le t<i\}|,\\[5pt]
c(\alpha)&=w_{f_k}-j-|\{t\,|\,w_t>t,w_t>j\}|+1.
\end{align*}
Then by \eqref{lmu}, we have
\begin{align}\label{hj}
\lambda_{r(\alpha)}+f_{r(\alpha)}-c(\alpha)+1
&=w_{f_k}-k+r(\alpha)-c(\alpha)+1\nonumber\\[5pt]
&=j-k+|\{t\,|\,w_t>t,w_t>j\}|+i-|\{t\,|\,w_t\le t<i\}|\nonumber\\[5pt]
&=j-|\{t\,|\,t< w_t\le j\}|-|\{t\,|\,w_t\le t<i\}|+i,
\end{align}
where, at the last step, we used the relation
\begin{align*}
k-|\{t\,|\,w_t>t,w_t>j\}|&=|\{t\,|\,w_t>t\}|-|\{t\,|\,w_t>t,w_t>j\}|\\[5pt]
&=|\{t\,|\,t< w_t\le j\}|.
\end{align*}
Since $w$ is 321-avoiding, it is easy to check that if there exists some integer $t$ such that $t<w_t\le j$, then $t<i$. Moreover, if $w_t\le t<i$, then  $w_t\le j$. Thus we have
\begin{align}\label{gf}
&j-|\{t\,|\,t< w_t\le j\}|-|\{t\,|\,w_t\le t<i\}|\nonumber\\[5pt]
&=j-(|\{t\,|\,t< w_t\le j,t<i\}|+|\{t\,|\,w_t\le t<i,w_t\le j\}|)\nonumber\\[5pt]
&=j-|\{t\,|\,t<i,w_t\le j\}|=|\{t\,|\,t\ge i,w_t\le j\}|\nonumber\\[5pt]
&=|\{(i,k)\in D(w)\,|\,k\leq j\}|\nonumber\\[5pt]
&=m_{ij}(w).
\end{align}
Combining \eqref{hj} and \eqref{gf} yields \eqref{rbr}.
This completes the proof. \qed

\subsection{Isobaric divided difference operator}

Before proving Theorem \ref{prop1} and Theorem \ref{prop2}, we recall some basic definitions.
Given a permutation $w=w_1w_2\cdots w_n\in S_n$,  the length   $\ell(w)$ of $w$ is equal to the number of inversions of $w$, namely,
\[\ell(w)=|\{(w_i,w_j)\,|\, 1\leq i<j\leq n, w_i>w_j\}|.\]
It is clear that $\ell(w)=|D(w)|$.
Let  $s_i$   denote the simple transposition interchanging  $i$ and $i+1$. Notice that $ws_i$ is the permutation obtained from $w$ by
swapping   $w_i$ and $w_{i+1}$. The  divided difference operator $\partial_i$ on the ring $\mathbb{Z}[\x]=\mathbb{Z}[x_1,x_2,\ldots, x_n]$ of polynomials with integer coefficients is defined by
\[\partial_i f(\x)=\frac{f(\x)
    -s_if(\x)}{x_i-x_{i+1}},\]
where $f(\x)\in \mathbb{Z}[\x]$ and $s_if(\x)$ is obtained from $f(\x)$ by interchanging $x_i$ and $x_{i+1}$.
One can then define  the isobaric divided difference operator $\pi_i$ as
\[\pi_i f(\x)=\partial_i (1-x_{i+1})f(\x).\]

The double Grothendieck polynomial $\mathfrak{G}_{w}(\x,\y)$ for $w\in S_n$ can be
 defined  as follows.
For the longest permutation $w_0=n \,(n-1)\cdots   1$, set
\[\mathfrak{G}_{w_0}(\x,\y)=\prod_{i+j\leq n}(x_i+y_j-x_iy_j).\]
 For $w\neq w_0$,  choose a simple transposition $s_i$   such that $\ell(ws_i)=\ell(w)+1$, and  let
 \begin{align}\label{GD}
 \mathfrak{G}_{w}(\x,\y)=\pi_i \mathfrak{G}_{ws_i}(\x,\y),
 \end{align}
where the operator  $\pi_i$ only acts  on the $\x$-variables.
Note that \eqref{GD} is independent of the choice of
the simple transposition $s_i$, since the operator
$\pi_i$ satisfies the Coxeter relations
 $\pi_i\pi_{i+1}\pi_i=\pi_{i+1}\pi_{i}\pi_{i+1}$ and
 $\pi_i\pi_{j}=\pi_j\pi_i$ for $|i-j|>1$.
If we set $y_i=0$ for $i\geq 1$, then $\mathfrak{G}_{w}(\x,\y)$ reduces to
the single Grothendieck polynomial $\mathfrak{G}_{w}(\x)$.

The double Schubert   polynomial $\mathfrak{S}_{w}(\x, \y)$ can be produced by a similar procedure \cite{LaSh-S,Ma}. Set
\[\mathfrak{S}_{w_0}(\x,\y)=\prod_{i+j\leq n}(x_i-y_j).\]
For $w\neq w_0$,  choose a simple transposition $s_i$   such that $\ell(ws_i)=\ell(w)+1$, and  let $\mathfrak{S}_{w}(\x,\y)=\partial_i \mathfrak{S}_{ws_i}(\x,\y)$.
By definition, it is easily seen that  $\mathfrak{S}_{w}(\x,\y)$
can be obtained from $\mathfrak{G}_{w}(\x,\y)$ by taking the lowest degree
homogeneous component and then replacing  $y_i$ by $-y_i$ for $i\geq 1$.
Analogously, putting $y_i=0$, $\mathfrak{S}_{w}(\x,\y)$ reduces to
the single Schubert  polynomial $\mathfrak{S}_{w}(\x)$. For combinatorial
constructions of Schubert polynomials, see for example \cite{As,AsSe,BeBi,BeSo,BiSt,FoRe,FK2,FoSt,LLS,WeYo,Wi}.

\subsection{Proof of Theorem \ref{prop1}}

The key idea is to show that, when $w$ is
 1432-avoiding, $G_w(\x,\y)$ is   compatible with
 the isobaric divided difference operator, which allows us to
 finish the proof by induction. Such an idea was first
 used by Wachs \cite{Wa} to establish the tableau formula for the Schubert polynomials
 of  2143-avoiding permutations.  Matsumura \cite{Matsumura-3} and
 Matsumura and  Sugimoto \cite{Matsumura-4} extended this idea
 to reprove the Knutson-Miller-Yong set-valued tableau formula for  the Grothendieck polynomials of
2143-avoiding permutations. Our technique  can be viewed as
 a generalization of that in  \cite{Matsumura-3,Matsumura-4} from
 Young diagrams to Rothe diagrams.

The longest permutation $w_0=n\cdots 2 1$ is  1432-avoiding. Since $D(w_0)$ is a staircase Young diagram with $n-i$ squares in row $i$, there is only one  tableau $T_0$ of shape $D(w_0)$ flagged  by $\f_0$, that is, every square  in the $i$-th row
of $T_0$ is filled with   $\{i\}$. For each square $(i,j)$ of $D(w_0)$, one has $m_{ij}(w_0)=j$.
 Thus,
\[
G_{w_0}(\x,\y)=(-1)^{|T_0|-\ell(w_0)}\prod_{i+j\leq n}(x_i\oplus y_j)=\prod_{i+j\leq n}(x_i\oplus y_j),
\]
which agrees with $\mathfrak{G}_{w_0}(\x,\y)$.

We now consider a 1432-avoiding permutation    $w\neq w_0$.
Let $r$ be the first ascent of $w$, that is,  the smallest position
such that $w_r<w_{r+1}$.
Lemma \ref{tp} claims  that
$ws_r$ is also   1432-avoiding.
Moreover, we will prove  that
\[G_w(\x,\y)=\pi_r G_{ws_r}(\x,\y),\]
see Theorem \ref{main-s}.
This allows us to give a proof of
 Theorem \ref{prop1} by  induction.

\begin{lemma}\label{tp}
Let $w\neq w_0$ be a 1432-avoiding permutation, and  $r$ be the first ascent of  $w$.
Then $ws_r$ is  a 1432-avoiding permutation.
\end{lemma}

\pf Write $w'=ws_r=w'_1w'_2\cdots w'_n$.
Suppose otherwise that $w'$ has a subsequence that is order isomorphic to $1432$.
Since $w$ is 1432-avoiding and $r$ is the first ascent, any subsequence of $w'$ that is order isomorphic to 1432 must be of the form  $w'_iw'_rw'_{r+1}w'_j$, where $i<r$ and $j>r+1$.
Since  $w'_i$ is the smallest element in this subsequence,
we have $w'_i<w'_{r+1}$. Noticing that  $w'_i=w_i$ and $w'_{r+1}=w_r$,
we see that $w_i<w_r$. However, since $r$ is the first ascent, we must have    $w_i>w_r$, leading to a contradiction.
This completes the proof. \qed

\begin{theorem}\label{main-s}
Let $w\neq w_0$ be a 1432-avoiding permutation, and $r$ be the first ascent of $w$. Then
\begin{align}\label{main-se}
G_{w}(\x,\y)=\pi_r G_{ws_r}(\x,\y).
\end{align}
\end{theorem}

In the rest of this subsection, we present a proof of Theorem \ref{main-s},
which can be outlined as follows.
We first define an equivalence relation  on the two sets $\SVRT(ws_r,\f_0)$ and $\SVRT(w,\f_0)$.
For an equivalence  class $C$ of $\SVRT(ws_r,\f_0)$, let
\begin{align}
G_{ws_r}(C; \x,\y)=(-1)^{\ell(ws_r)}\sum_{T\in C} \prod_{(i,j)\in D(w)}\prod_{t\in T({i,j})}(-x_t\oplus y_{m_{ij}(w)+i-t})\label{class-g}
\end{align}
denote the polynomial generated by the Rothe tableaux in   $C$.
In Theorem \ref{mm-1}, we deduce a formula for $G_{ws_r}(C; \x,\y)$.
Similarly,  write $G_{w}(C'; \x,\y)$ for the polynomial generated by the Rothe tableaux in an equivalence class $C'$ of $\SVRT(w, \f_0)$. We also obtain an expression for  $G_{w}(C'; \x,\y)$, see Theorem \ref{mm-2}. Finally, we establish    a bijection $\Phi$ between the set of equivalence classes of $\SVRT(ws_r,\f_0)$ and
the set of equivalence classes of $\SVRT(w, \f_0)$. The formulas given in Theorems \ref{mm-1} and
 \ref{mm-2} allow us to conclude  that for any equivalence class $C$ of $\SVRT(ws_r,\f_0)$,
\[\pi_r G_{ws_r}(C; \x,\y)=G_w(\Phi(C); \x,\y).\]
This  leads to a proof of Theorem  \ref{main-s}.

Unless otherwise stated, we  always assume that $w\neq w_0$ is a 1432-avoiding permutation,  and  that $r$ is the
first ascent of $w$. For  $T\in\SVRT(ws_r, \f_0)$, let
\begin{align*}
E(T)=\{B\in D(ws_r)\,|\,\{r,r+1\}\cap T(B)\neq\emptyset\},
\end{align*}
that is, $E(T)$ is the subset  of squares of $T$ containing at least one of  $r$ and $r+1$. It should be noted that  the definition of $E(T)$ for $T\in\SVRT(ws_r, \f_0)$ depends only on $r$, which has nothing to do with the first ascent of $ws_r$.

\begin{defi}
 Given two Rothe tableaux $T,T'\in\SVRT(ws_r, \f_0)$, we say that $T$ is equivalent to $T'$, denoted $T\sim T'$,
 if $E(T)=E(T')$ and  for every square $B\in D(ws_r)$,
\[T(B)\setminus \{r,r+1\}=T'(B)\setminus \{r,r+1\}.\]
\end{defi}
The equivalence relation on the set $\SVRT(w, \f_0)$ is defined in the same manner. Let $\SVRT(ws_r, \f_0)/\hspace{-.1cm}\sim$ and $\SVRT(w, \f_0)/\hspace{-.1cm}\sim$  denote the sets of equivalence classes of $\SVRT(ws_r, \f_0)$ and $\SVRT(w, \f_0)$, respectively.

Given a Rothe tableau $T$ in $\SVRT(ws_r, \f_0)$ or $\SVRT(w, \f_0)$, since the columns of $T$ are strictly increasing, each column of $T$ contains at most two squares  in $E(T)$.
Let $P(T)$ be the subset of $E(T)$ such that a square $B\in E(T)$ belongs to
$P(T)$ if the   column of $T$ containing $B$ has only one square (i.e., $B$) in $E(T)$.
Let $Q(T)=E(T)\setminus P(T)$, namely, the subset of $E(T)$ such that
a square $B\in E(T)$ belongs to
$Q(T)$ if the   column containing $B$ has exactly two squares in $E(T)$.
Evidently, $T\sim T'$ if and only if
\[P(T)=P(T')\quad \text{and}\quad Q(T)=Q(T').\]

Let $C$ be an equivalence class of $\SVRT(ws_r, \f_0)$ or  $\SVRT(w, \f_0)$, and
let $T$ be any given Rothe tableau in $C$. For $i\ge1$, let $P(T,i)$ be the set of squares of $P(T)$ in row $i$, and  let $b_i(T)=|P(T,i)|$. Clearly, $P(T,i)$ is empty unless $i\geq r$.
Moreover,  let
\begin{align}\label{ellr}
\ell_i(T)=m_{ip_i}(ws_r)+i-r-1,
\end{align}
where $(i,p_i)$ is the leftmost square in  $P(T,i)$. To state the formula for  $G_{ws_r}(C; \x,\y)$ or $G_{w}(C; \x,\y)$, we need to define a polynomial $h(C,i;\x,\y)$.
Set $h(C,i;\x,\y)=1$ if $b_i(T)=0$, and for $b_i(T)\geq 1$, let
\begin{align}
h(C,i;\x,\y)=&\sum_{k=0}^{b_i(T)}\, \prod_{j=1}^{k} (-x_{r+1}\oplus y_{\ell_i(T)+j-1})
\prod_{j=k+1}^{b_i(T)} (-x_{r}\oplus y_{\ell_i(T)+j})\nonumber\\[5pt]
&+\sum_{k=1}^{b_i(T)}\, \prod_{j=1}^k (-x_{r+1}\oplus y_{\ell_i(T)+j-1})
\prod_{j=k}^{b_i(T)} (-x_{r}\oplus y_{\ell_i(T)+j}).\label{hc}
\end{align}
Note that $h(C,i;\x,\y)$ is independent of the choice of the Rothe tableau $T$ in $C$.
As will be seen in the proof  of Theorem \ref{mm-1}, $h(C,i;\x,\y)$ records  the contributions of the integers  $r$ and $r+1$ in the squares of $P(T,i)$ ($i\ge r+1$)
summed over all the Rothe tableaux in $C$.

\begin{theorem}\label{mm-1}
Let $w\neq w_0$ be a 1432-avoiding permutation, and $r$ be the first ascent of $w$.
Assume that $C\in\SVRT(ws_r, \f_0)/\hspace{-.1cm}\sim$ and  $T$ is any given Rothe tableau in $C$.
Then,
\begin{align}
G_{ws_r}(C; \x,\y)=&(-1)^{\ell(ws_r)} \left(\prod_{(i,j)\in D(ws_r)}\prod_{t\in T(i,j)\atop t\neq r,r+1}(-x_t\oplus y_{m_{ij}(ws_r)+i-t})\right)\nonumber \\[5pt]
&\ \ \ \ \ \cdot\left(\prod_{j=1}^{b_r(T)}(-x_r\oplus y_{\ell_r(T)+j})\right)\cdot H_C(\x,\y)\cdot  J_C(\x,\y).\label{gggs}
\end{align}
In the above expression \eqref{gggs},
\[H_C(\x,\y)=\prod_{ i>r+1}h(C,i;\x,\y),\]
and
\[J_C(\x,\y)=\prod_{(i,j)\in Q^+(T)}(-x_r\oplus y_{m_{ij}(ws_r)+i-r})(-x_{r+1}\oplus y_{m_{ij}(ws_r)+i-r}),\]
where $Q^{+} (T)$
 denotes the subset of $Q(T)$
consisting of the squares containing $r$.
\end{theorem}

Although the formula for $G_{ws_r}(C,\x,\y)$ in \eqref{gggs} looks  a bit complicated, it will be clear from the proof that each factor in \eqref{gggs} appears naturally.
In fact, for two different Rothe tableaux $T,T'\in C$, $T$ and $T'$ can be possibly different only in the squares of $P(T)$. Thus the first factor is the contribution of the integers other than $r$ and $r+1$. We shall show that the second factor is the contribution of $r$ in $P(T,r)$, $J_C(\x,\y)$ is the contribution of $r$ and $r+1$ in $Q(T)$, and $H_C(\x,\y)$ is the contribution of $r$ and $r+1$ in $P(T,i)$ with $i\ge r+1$ summed over all the Rothe tableaux in $C$.

To prove Theorem \ref{mm-1}, we need two lemmas concerning
 the configuration of the squares in the sets $P(T)$ and $Q(T)$.

\begin{lemma}\label{right}
Let $w\neq w_0$ be a 1432-avoiding permutation.
Assume that  $T$ is a Rothe tableau in $\SVRT(w, \f_0)$ or $\SVRT(ws_r, \f_0)$, and  $(i,j)\in P(T)$.   Then there do not exist two squares $(i,k),(h,k)\in Q(T)$ such that $k>j$ and $h<i$.
\end{lemma}

\begin{figure}[h]
\setlength{\unitlength}{0.5mm}
\begin{center}
\begin{picture}(100,60)
\put(50,0){\line(1,0){10}}\put(50,10){\line(1,0){10}}
\put(50,0){\line(0,1){10}}\put(60,0){\line(0,1){10}}


\put(80,20){\line(1,0){10}}\put(80,30){\line(1,0){10}}
\put(80,20){\line(0,1){10}}\put(90,20){\line(0,1){10}}

\put(80,0){\line(1,0){10}}\put(80,10){\line(1,0){10}}
\put(80,0){\line(0,1){10}}\put(90,0){\line(0,1){10}}

\qbezier[50](20,5)(33,5)(46,5)\qbezier[50](20,25)(33,25)(46,25)
\put(10,4){$i$}\put(10,23.5){$h$}

\qbezier[25](55,34)(55,40)(55,46)\qbezier[25](85,34)(85,40)(85,46)
\put(54,52){$j$}\put(83.5,51.5){$k$}

\end{picture}
\end{center}
\caption{An illustration for the proof of Lemma \ref{right}.}
\label{a}
\end{figure}
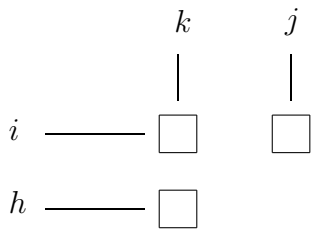

\pf We   only give a proof for the case when $T\in\SVRT(w, \f_0)$. The same argument applies to the case when $T\in \SVRT(ws_r, \f_0)$.
Suppose to the contrary that there exist two squares $(i,k),(h,k)\in Q(T)$ such that $k>j$ and $h<i$, see Figure \ref{a} for  an illustration.

Since   $(i,j),(h,k)\in D(w)$, it follows that $w_h>j$ and $j$ appears after $w_h$ in $w$.
This implies that $(h,j)\in D(w)$.
Keep in mind that each of the sets $T(i,j)$, $T(i,k)$ and $T(h,k)$ contains at least one of the integers  $r$ and $r+1$. Since the rows of $T$ are weakly decreasing and the   columns of $T$ are strictly increasing, we see that $r\in T(h,k)$, $r+1\in T(i,j)$. This forces  that $T(h,j)=\{r\}$, and hence    $(i,j)\in Q(T)$,
which contradicts   the assumption that $(i,j)\in P(T)$.
\qed

It should be noted that Lemma \ref{right} is valid for any permutation since the pattern avoidance condition is not required in  the proof.

\begin{lemma}\label{left}
Let $w\neq w_0$ be a 1432-avoiding permutation.
Assume that $T$ is a Rothe tableau in $\SVRT(w, \f_0)$  or  $\SVRT(ws_r, \f_0)$, and  $(i,j)\in P(T)$. If $i>r$, then
there do not exist two squares $(i,k),(h,k)\in Q(T)$
such that $h>i$ and $k<j$.
\end{lemma}

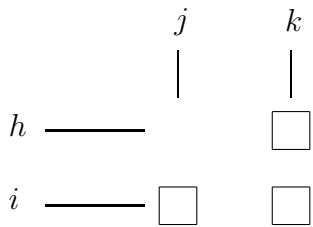
\begin{figure}[h]
\setlength{\unitlength}{0.5mm}
\begin{center}
\begin{picture}(100,60)
\put(50,0){\line(1,0){10}}\put(50,10){\line(1,0){10}}
\put(50,0){\line(0,1){10}}\put(60,0){\line(0,1){10}}

\put(50,20){\line(1,0){10}}\put(50,30){\line(1,0){10}}
\put(50,20){\line(0,1){10}}\put(60,20){\line(0,1){10}}

\put(80,20){\line(1,0){10}}\put(80,30){\line(1,0){10}}
\put(80,20){\line(0,1){10}}\put(90,20){\line(0,1){10}}

\qbezier[50](20,5)(33,5)(46,5)\qbezier[50](20,25)(33,25)(46,25)
\put(10,4){$h$}\put(10,23.5){$i$}

\qbezier[25](55,34)(55,40)(55,46)\qbezier[25](85,34)(85,40)(85,46)
\put(54,52){$k$}\put(83.5,53){$j$}

\end{picture}
\end{center}
\caption{An illustration for the proof of Lemma \ref{left}.}
\label{b}
\end{figure}

\pf We only give  a proof for   $T\in \SVRT(w, \f_0)$, and
the arguments for  $T\in \SVRT(ws_r, \f_0)$ can be carried out in the same manner. Suppose  otherwise that there exist two   squares $(i,k)$ and $(h,k)$ in $Q(T)$ where $i<h$ and $j>k$, as illustrated in Figure \ref{b}.
Notice  that both the sets $T(i,j)$ and $T(i,k)$ contain $r$, while the set $T(h,k)$ contains $r+1$.
We have the following two claims.

\noindent
Claim 1: $w_s<k$ for  any $i<s<h$.
 Suppose otherwise  that there exists some $i<s<h$ such that $w_s>k$.
 Then the square $(s,k)$ belongs to $D(w)$. Since $r\in T(i,k)$
  and $r+1\in T(h,k)$, it follows that
 $r<\min T(s,k)<r+1$, leading to a contradiction.

\noindent
Claim 2:   $k<w_h\leq j$. Since $(h,k)$ is a square in $D(w)$, it is clear that
$k<w_h$.
Suppose otherwise that $w_h>j$. It
 follows from Claim 1 that $j$ must appear in $w$ after the position $h$.
 This implies that  $(h,j)\in D(w)$.
 Since $r\in T(i,j)$ and $r+1\in T(h,k)$,
we must have $T(h,j)=\{r+1\}$. This implies  that $(i,j)\in Q(T)$,
 contradicting the assumption that $(i,j)\in P(T)$.

By Claim 2 and the fact that $w_i>j$, we see that   $w_iw_h k$ forms a decreasing subsequence of  $w$.
Since $w$ is 1432-avoiding, we have $w_t> k$ for any $1\leq t<i$.
Thus, for any $1\leq t< i$, the square $(t,k)$ belongs to $D(w)$. Keep in mind that each integer in row $i$ of $T$ cannot exceed $i$ and the columns of $T$ are strictly increasing. So we have
 $T(t,k)=\{t\}$ for $1\leq t\leq i$. In particular, we have $T(i,k)=\{i\}$.
Since $r\in T(i,k)$, we must have $i=r$, contradicting  the
assumption that $i>r$. This completes the proof. \qed

Based on Lemmas \ref{right} and \ref{left},
we can now give a proof of
Theorem \ref{mm-1}.

\noindent
\textit{Proof of Theorem \ref{mm-1}.}
Assume that $T'\in \SVRT(ws_r,\f_0)$ is a Rothe tableau in the equivalence class  $C$. Then $T'$ differs from $T$  only possibly in
the squares  of $P(T)$. Note that if $P(T,i)$ is nonempty, then we must have $i\geq r$.
Moreover, since the integers appearing in $r$-th row of  $T'$
cannot exceed $r$, it follows that for any square $B\in P(T,r)$, $T'(B)$ does not contain $r+1$.
Thus, for $B\in P(T,r)$, $r\in T(B)=T'(B)$ and $r+1\notin T(B)=T'(B)$.

\begin{figure}[h]
\setlength{\unitlength}{0.5mm}
\begin{center}
\begin{picture}(270,60)
\put(0,0){\line(1,0){20}}\put(20,0){\line(0,1){20}}
\put(20,20){\line(1,0){20}}\put(40,20){\line(0,1){10}}
\put(40,30){\line(1,0){20}}\put(60,30){\line(0,1){10}}
\put(60,40){\line(1,0){30}}\put(90,40){\line(0,1){10}}
\put(90,50){\line(-1,0){90}}\put(0,0){\line(0,1){50}}

\put(30,0){\line(1,0){10}}\put(30,0){\line(0,1){10}}
\put(40,0){\line(0,1){10}}\put(30,10){\line(1,0){10}}

\put(50,0){\line(1,0){10}}\put(50,10){\line(1,0){10}}
\put(50,0){\line(0,1){10}}\put(60,0){\line(0,1){10}}

\put(70,0){\line(1,0){20}}\put(70,10){\line(1,0){20}}
\put(70,0){\line(0,1){10}}\put(80,0){\line(0,1){10}}
\put(90,0){\line(0,1){10}}

\put(100,0){\line(1,0){10}}\put(100,10){\line(1,0){10}}
\put(100,0){\line(0,1){10}}\put(110,0){\line(0,1){10}}

\qbezier[20](-9,15)(-5.5,15)(-2,15)\put(-17,13){ $r$}

\put(95,45){\circle*{3}}\put(65,35){\circle*{3}}
\put(45,25){\circle*{3}}\put(115,5){\circle*{3}}
\put(25,15){\circle*{3}}


\put(160,0){\line(1,0){20}}\put(180,0){\line(0,1){20}}
\put(180,20){\line(1,0){20}}\put(200,20){\line(0,1){10}}
\put(200,30){\line(1,0){20}}\put(220,30){\line(0,1){10}}
\put(220,40){\line(1,0){30}}\put(250,40){\line(0,1){10}}
\put(250,50){\line(-1,0){90}}\put(160,0){\line(0,1){50}}

\put(180,10){\line(1,0){10}}

\put(190,10){\line(1,0){10}}\put(190,20){\line(1,0){10}}
\put(190,10){\line(0,1){10}}\put(200,10){\line(0,1){10}}

\put(210,10){\line(1,0){10}}\put(210,20){\line(1,0){10}}
\put(210,10){\line(0,1){10}}\put(220,10){\line(0,1){10}}

\put(230,10){\line(1,0){20}}\put(230,20){\line(1,0){20}}
\put(230,10){\line(0,1){10}}\put(240,10){\line(0,1){10}}
\put(250,10){\line(0,1){10}}

\put(260,10){\line(1,0){10}}\put(260,20){\line(1,0){10}}
\put(260,10){\line(0,1){10}}\put(270,10){\line(0,1){10}}

\qbezier[20](151,15)(155.5,15)(158,15)\put(143,13){ $r$}

\put(255,45){\circle*{3}}\put(225,35){\circle*{3}}
\put(205,25){\circle*{3}}\put(275,15){\circle*{3}}
\put(185,5){\circle*{3}}

\put(182,12){$\heartsuit$}

\end{picture}
\end{center}
\caption{The first $r+1$ rows of $D(w)$ and $D(ws_r)$.}
\label{zao}
\end{figure}
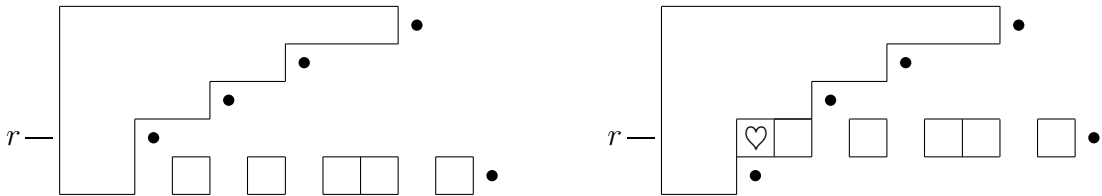

Before we proceed, we give an illustration of the configuration of the
squares in the first $r+1$ rows of $D(w)$ and $D(ws_r)$, which will be helpful to
analyze the contributions of the integer $r$  in the squares of $P(T,r)$.
Notice that    $D(w)$ is obtained from  $D(ws_r)$ by deleting the square $(r,w_r)$ and then moving each square in row $r$, that lies  to the right of $(r,w_r)$, down to row $r+1$. Since $r$ is the first ascent of $w$, the first $r+1$ rows of  $D(w)$
and $D(ws_r)$ are as depicted  in Figure \ref{zao}, where the square $(r,w_r)$ of $D(ws_r)$ is signified  by a symbol $\heartsuit$.

Obviously, the first $w_r-1$ squares in the $r$-th row (respectively, ($r+1$)-th row) of  $T$ are filled with the set  $\{r\}$ (respectively, $\{r+1\}$). This implies that each square in the ($r+1$)-th row of $D(ws_r)$ belongs to $Q(T)$ and
the set  $P(T, r+1)$ is empty.
Therefore,  the contribution of the $r$'s in squares of $P(T,r)$ to $G_{ws_r}(C; \x,\y)$ is
\begin{align}\label{ft}
\prod_{(r,j)\in P(T,r)}(-x_r\oplus y_{m_{rj}(ws_r)}).
\end{align}

On the other hand, the contribution of the $r$'s and $(r+1)$'s in squares of $P(T,i)$ for $i>r+1$ to $G_{ws_r}(C; \x,\y)$ is
\begin{align}
F_C(\x,\y)=&\sum_{T'\in C}\ \prod_{i> r+1}\prod_{(i,j)\in P^+(T',i)}(-x_r\oplus y_{m_{ij}(ws_r)+i-r})\nonumber\\[5pt]
&\ \ \cdot\prod_{(i,j)\in P^-(T',i)} \ (-x_{r+1}\oplus y_{m_{ij}(ws_r)+i-r-1}),\label{tie}
\end{align}
where $P^+(T',i)$ (respectively, $P^-(T',i)$) denotes  the subset of  $P(T',i)$ consisting of squares containing $r$ (respectively, $r+1$). Moreover, the contribution of the $r$'s and $(r+1)$'s in squares of $Q(T)$ to $G_{ws_r}(C; \x,\y)$ is
\begin{align}\label{Qcontri}
R_C(\x,\y)=\prod_{(i,j)\in Q^+(T)}(-x_r\oplus y_{m_{ij}(ws_r)+i-r})\prod_{(i,j)\in Q^-(T)}(-x_{r+1}\oplus y_{m_{ij}(ws_r)+i-r-1}),
\end{align}
where $Q^{+}(T)$ (respectively, $Q^-(T)$)
 denotes the subset of $Q(T)$
consisting of the squares containing $r$ (respectively, $r+1$).
Consequently, we obtain that
\begin{align}
G_{ws_r}(C; \x,\y)=&(-1)^{\ell(ws_r)}\left( \prod_{(i,j)\in D(ws_r)}\prod_{t\in T(i,j)\atop t\neq r,r+1}(-x_t\oplus y_{m_{ij}(ws_r)+i-t})\right)\nonumber \\[5pt]
&\ \ \cdot\prod_{(r,j)\in P(T,r)}(-x_r\oplus y_{m_{rj}(ws_r)})\cdot F_C(\x,\y)\cdot R_C(\x,\y).\label{kong}
\end{align}

Comparing \eqref{kong} with \eqref{gggs}, in order to complete the proof, we need to show that
\begin{align}\label{fit}
\prod_{(r,j)\in P(T,r)}(-x_r\oplus y_{m_{rj}(ws_r)})&=\prod_{j=1}^{b_r(T)}(-x_r\oplus y_{\ell_r(T)+j}),\\[5pt]
F_C(\x,\y)&=\prod_{ i>r+1}h(C,i;\x,\y)=H_C(\x,\y),\label{f1}\\[5pt]
R_C(\x,\y)&=J_C(\x,\y).\label{f2fs}
\end{align}

Let us first prove  \eqref{fit}. To this end, we show that  if there are two squares $(r,j_1)$ and $(r,j_2)$ in $P(T,r)$ with $j_1<j_2$  and there exists a square $(r,j)\in D(ws_r)$ for some $j_1<j<j_2$, then $(r,j)\in P(T,r)$.
It suffices to prove the following claim.

\noindent
{\bf Claim.} For $(r,j)\in P(T, r)$,
there do not exist squares $(r,k),(h,k)\in Q(T)$
such  that $h>r+1$ and $k<j$.

To verify this claim, we construct a Rothe tableau $\overline{T}$ from $T$ such that $\overline{T}\in \SVRT(w, \f_0)$.
Let $R$ be the set of squares of $D(ws_r)$ in row $r$ that are strictly to the  right of  $(r,w_r)$.
Define $\overline{T}$ to be the tableau obtained from $T$ by  deleting the square $(r,w_r)$ together with $T(r,w_r)$, and then moving each square $B$ in $R$, together with $T(B)$, down to row $r+1$. By  construction, it is easy to check that $\overline{T}\in\SVRT(w, \f_0)$.
Note that  $(r,j)\in P(T)$ if and  only if
  $(r+1,j)\in P(\overline{T})$. Applying Lemma \ref{left}
 to $\overline{T}$, we see that if $(r+1,j)\in P(\overline{T})$, then there do not exist squares $(r,k),(h,k)\in Q(\overline{T})$
with $h>r+1$ and $k<j$.
Since $Q(T)=Q(\overline{T})$, we conclude  the claim.

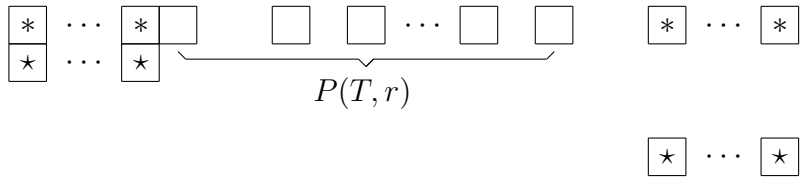
\begin{figure}[h]
\setlength{\unitlength}{0.5mm}
\begin{center}
\begin{picture}(185,45)

\put(-10,25){\line(1,0){10}}\put(-10,25){\line(0,1){10}}
\put(0,25){\line(0,1){10}}\put(-10,35){\line(1,0){10}}

\put(-10,15){\line(1,0){10}}\put(-10,15){\line(0,1){10}}
\put(0,15){\line(0,1){10}}\put(-10,25){\line(1,0){10}}

\put(20,25){\line(1,0){10}}\put(20,25){\line(0,1){10}}
\put(30,25){\line(0,1){10}}\put(20,35){\line(1,0){10}}

 \put(20,15){\line(1,0){10}}\put(20,15){\line(0,1){10}}
\put(30,15){\line(0,1){10}}\put(20,25){\line(1,0){10}}

\put(30,25){\line(1,0){10}}\put(30,25){\line(0,1){10}}
\put(40,25){\line(0,1){10}}\put(30,35){\line(1,0){10}}

\put(60,25){\line(1,0){10}}\put(60,25){\line(0,1){10}}
\put(70,25){\line(0,1){10}}\put(60,35){\line(1,0){10}}

\put(80,25){\line(1,0){10}}\put(80,25){\line(0,1){10}}
\put(90,25){\line(0,1){10}}\put(80,35){\line(1,0){10}}

\put(110,25){\line(1,0){10}}\put(110,25){\line(0,1){10}}
\put(120,25){\line(0,1){10}}\put(110,35){\line(1,0){10}}

\put(130,25){\line(1,0){10}}\put(130,25){\line(0,1){10}}
\put(140,25){\line(0,1){10}}\put(130,35){\line(1,0){10}}

\put(35,23){\line(1,-1){2}}\put(37,21){\line(1,0){46}}
\put(83,21){\line(1,-1){2}}\put(85,19){\line(1,1){2}}
\put(85,19){\line(1,1){2}}\put(87,21){\line(1,0){46}}
\put(133,21){\line(1,1){2}}
\put(71,10){$P(T,r)$}

\put(160,25){\line(1,0){10}}\put(160,25){\line(0,1){10}}
\put(170,25){\line(0,1){10}}\put(160,35){\line(1,0){10}}

 \put(160,-10){\line(1,0){10}}\put(160,-10){\line(0,1){10}}
\put(170,-10){\line(0,1){10}}\put(160,0){\line(1,0){10}}

\put(190,25){\line(1,0){10}}\put(190,25){\line(0,1){10}}
\put(200,25){\line(0,1){10}}\put(190,35){\line(1,0){10}}

\put(190,-10){\line(1,0){10}}\put(190,-10){\line(0,1){10}}
\put(200,-10){\line(0,1){10}}\put(190,0){\line(1,0){10}}

\put(5,28){$\cdots$}\put(5,18){$\cdots$}
\put(175,28){$\cdots$}\put(175,-7){$\cdots$}
\put(95,28){$\cdots$}

\put(-7,28){$\ast$}\put(-7,18){$\star$}
\put(23,28){$\ast$}\put(23,18){$\star$}

\put(163,-7){$\star$}\put(163,28){$\ast$}
\put(193,-7){$\star$}\put(193,28){$\ast$}

\end{picture}
\end{center}
\caption{An illustration of the squares in $P(T,r)$.}
\label{anm11}
\end{figure}

By the above claim, the configuration  of the squares of $P(T)$ and $Q(T)$ in the $r$-th row of $D(ws_r)$ is as illustrated in
Figure \ref{anm11}, where the squares in $Q^+(T)$
(respectively, $Q^-(T)$) are marked with a $\ast$ (respectively, $\star$). In view of the definition $m_{i,j}(ws_r)$ in \eqref{statistic} as well as the definition $\ell_r(T)$ in \eqref{ellr}, we see that \eqref{fit} holds.

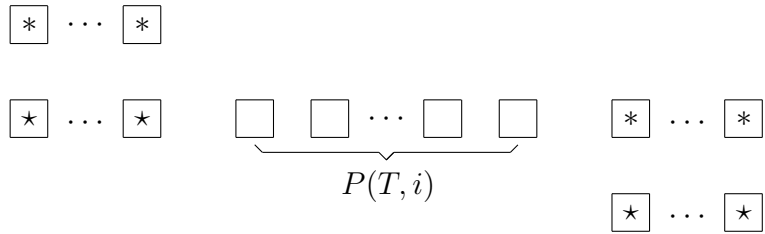
\begin{figure}[h]
\setlength{\unitlength}{0.5mm}
\begin{center}
\begin{picture}(200,70)
\put(0,25){\line(1,0){10}}\put(0,25){\line(0,1){10}}
\put(10,25){\line(0,1){10}}\put(0,35){\line(1,0){10}}

 \put(0,50){\line(1,0){10}}\put(0,50){\line(0,1){10}}
\put(10,50){\line(0,1){10}}\put(0,60){\line(1,0){10}}

\put(30,25){\line(1,0){10}}\put(30,25){\line(0,1){10}}
\put(40,25){\line(0,1){10}}\put(30,35){\line(1,0){10}}

 \put(30,50){\line(1,0){10}}\put(30,50){\line(0,1){10}}
\put(40,50){\line(0,1){10}}\put(30,60){\line(1,0){10}}

\put(60,25){\line(1,0){10}}\put(60,25){\line(0,1){10}}
\put(70,25){\line(0,1){10}}\put(60,35){\line(1,0){10}}

\put(80,25){\line(1,0){10}}\put(80,25){\line(0,1){10}}
\put(90,25){\line(0,1){10}}\put(80,35){\line(1,0){10}}

\put(110,25){\line(1,0){10}}\put(110,25){\line(0,1){10}}
\put(120,25){\line(0,1){10}}\put(110,35){\line(1,0){10}}

\put(130,25){\line(1,0){10}}\put(130,25){\line(0,1){10}}
\put(140,25){\line(0,1){10}}\put(130,35){\line(1,0){10}}

\put(65,23){\line(1,-1){2}}\put(67,21){\line(1,0){31}}
\put(98,21){\line(1,-1){2}}\put(100,19){\line(1,1){2}}
\put(100,19){\line(1,1){2}}\put(102,21){\line(1,0){31}}
\put(133,21){\line(1,1){2}}
\put(88,10){$P(T,i)$}

\put(160,25){\line(1,0){10}}\put(160,25){\line(0,1){10}}
\put(170,25){\line(0,1){10}}\put(160,35){\line(1,0){10}}

 \put(160,0){\line(1,0){10}}\put(160,0){\line(0,1){10}}
\put(170,0){\line(0,1){10}}\put(160,10){\line(1,0){10}}

\put(190,25){\line(1,0){10}}\put(190,25){\line(0,1){10}}
\put(200,25){\line(0,1){10}}\put(190,35){\line(1,0){10}}

 \put(190,0){\line(1,0){10}}\put(190,0){\line(0,1){10}}
\put(200,0){\line(0,1){10}}\put(190,10){\line(1,0){10}}

\put(15,53){$\cdots$}\put(15,27){$\cdots$}
\put(175,27){$\cdots$}\put(175,2){$\cdots$}
\put(95,28){$\cdots$}

\put(3,28){$\star$}\put(3,53){$\ast$}
\put(33,28){$\star$}\put(33,53){$\ast$}

\put(163,3){$\star$}\put(163,28){$\ast$}
\put(193,3){$\star$}\put(193,28){$\ast$}

\end{picture}
\end{center}
\caption{The configuration  of the squares in $P(T,i)$ with $i>r+1$.}
\label{anm}
\end{figure}

We next prove    \eqref{f1}.
For $i>r$, by Lemma \ref{right} and Lemma \ref{left}, the configuration
 of the squares
of $P(T)$ and $Q(T)$ must be as illustrated as in Figure \ref{anm}. In particular, every square in row $i$ of $D(ws_r)$ that lies between the leftmost square and the rightmost square of $P(T,i)$ must  belong to $P(T,i)$.
Assume that $T'$ is
  a Rothe tableau in $C$. Keep in mind that $P(T',i)=P(T,i)$. Then we have the following two
cases.

\noindent
Case 1:
In $T'$, the first  $k$ $(0\leq k\leq b_i(T))$ squares in $P(T,i)$ contain $r+1$, and the remaining  $b_i(T)-k$
squares in $P(T,i)$ contain $r$.
In this case, running over   the Rothe tableaux $T'$  in $C$, the integers $r$ and
$r+1$ in $P(T,i)$ contribute
\begin{align}\label{QF-1}
\sum_{k=0}^{b_i(T)}\prod_{j=1}^{k} (-x_{r+1}\oplus y_{\ell_i(T)+j-1})
\prod_{j=k+1}^{b_i(T)} (-x_{r}\oplus y_{\ell_i(T)+j}).
\end{align}

\noindent
Case 2: In  $T'$, the first  $k-1$ $(1\leq k\leq b_i(T))$ squares in $P(T,i)$ contain $r+1$,   the   $k$-th square contains both $r$ and $r+1$, and the remaining   $b_i(T)-k$ squares in
$P(T,i)$ contain $r$.
In this case, running  over  the Rothe tableaux $T'$  in $C$, the integers $r$ and
$r+1$ in $P(T,i)$ contribute
\begin{align}\label{QF-2}
\sum_{k=1}^{b_i(T)}\prod_{j=1}^{k} (-x_{r+1}\oplus y_{\ell_i(T)+j-1})
\prod_{j=k}^{b_i(T)} (-x_{r}\oplus y_{\ell_i(T)+j}).
\end{align}
Combining \eqref{QF-1} and \eqref{QF-2}, we see that
Case 1 and Case 2 together contribute the factor $h(C,i;\x,\y)$ as defined in \eqref{hc} to the summation $F_C(\x,\y)$ in \eqref{tie}.
Running over the row indices $i$ with $i>r+1$ yields  \eqref{f1}.

Finally, we verify \eqref{f2fs}. For each $(i,j)\in Q^+(T)$, we use $(i',j)$ to denote the square in $Q^-(T)$ that lie in the same column as $(i,j)$.
Then we have
\begin{align}\label{HD-1}
R_C(\x,\y)=\prod_{(i,j)\in Q^+(T)}(-x_r\oplus y_{m_{ij}(ws_r)+i-r})
(-x_{r+1}\oplus y_{m_{i'j}(ws_r)+i'-r-1}).
\end{align}
Write $w'=ws_r$. We assert that $w'_t<j$ for $i<t<i'$. Suppose
otherwise that $w'_t>j$. Since the square $(i',j)\in D(ws_r)$, we see that $(t,j)\in D(ws_r)$.
Thus we have $r<\min T(t,j)\leq \max T(t,j)<r+1$,
leading to a contradiction. This verifies the assertion.
By the definition of $m_{ij}(w)$ in \eqref{statistic}, it is easy to see that \[m_{ij}(w')=|\{(i,k)\in D(w')\,|\,k\le j\}|=|\{t>i\,|\,w'_t\le j\}|.\]
Therefore, by the above assertion, we obtain
\[m_{ij}(ws_r)=m_{i'j}(ws_r)+i'-i-1,\]
and so that
\begin{align}\label{HD-2}
m_{i'j}(ws_r)+i'-r-1=m_{ij}(ws_r)+i-r.
\end{align}
Putting \eqref{HD-2} into \eqref{HD-1}, we arrive at the equality in \eqref{f2fs}.
This completes the proof. \qed

The following theorem provides  a formula for the polynomial $G_{w}(C'; \x,\y)$ for an equivalence class $C'\in\SVRT(w, \f_0)/\hspace{-.1cm}\sim$.

\begin{theorem}\label{mm-2}
Let $w\neq w_0$ be a 1432-avoiding permutation, and $r$ be the first ascent of $w$.
Assume that $C'\in\SVRT(w, \f_0)/\hspace{-.1cm}\sim$ and  $T'$ is any given Rothe tableau in $C'$.
Then,
\begin{align}
G_{w}(C'; \x,\y)=&(-1)^{\ell(w)} \left(\prod_{(i,j)\in D(w)}\prod_{t\in T'(i,j)\atop t\neq r,r+1}(-x_t\oplus y_{m_{ij}(w)+i-t})\right)\nonumber \\[5pt]
&\ \cdot h(C',r+1;\x,\y)\cdot H_{C'}(\x,\y)\cdot J_{C'}(\x,\y),\label{gggg}
\end{align}
where
\[{H}_{C'}(\x,\y)=\, \prod_{ i>r+1}{h}(C',i;\x,\y),\]
and
\[J_{C'}(\x,\y)=\prod_{(i,j)\in Q^+(T')}(-x_r\oplus y_{m_{ij}(w)+i-r})(-x_{r+1}\oplus y_{m_{ij}(w)+i-r}).\]
\end{theorem}

\noindent
{\it Sketch of the proof.} The proof is nearly the same as the arguments for
Theorem \ref{mm-1}. The only difference is to notice that
$P(T',r)$ is empty  and that the squares in  $P(T',r+1)$ contributes
the factor $h(C',r+1;\x,\y)$. \qed

To finish the proof of Theorem \ref{main-s}, we need a one-to-one correspondence  $\Phi$
between  the two sets of equivalence classes:
\[
\Phi:\ \SVRT(ws_r,\f_0)/\hspace{-.1cm}\sim\ \longrightarrow\ \SVRT(w,\f_0)/\hspace{-.1cm}\sim.
\]

\noindent{\textbf{Construction of the bijection $\Phi$}:}
Assume that  $C\in\SVRT(ws_r, \f_0)/\hspace{-.1cm}\sim$ and  $T\in C$. Let $T'=\overline{T}\in\SVRT(w, \f_0)$ be the Rothe tableau as constructed in the proof of \eqref{fit}. That is,  $T'$ is the Rothe tableau obtained from $T$ by  deleting the square $(r,w_r)$ together with $T(r,w_r)$, and then moving each square $B$ in $R$, together with $T(B)$, down to row $r+1$, where $R$ is the set of squares of $D(ws_r)$ in row $r$ that are strictly to the  right of  $(r,w_r)$.
Let $C'\in\SVRT(w, \f_0)/\hspace{-.1cm}\sim$
be the equivalence class containing  $T'$. It is clear that
$C'$ is independent of the choice of $T$. Set $\Phi(C)=C'$.

The inverse  of $\Phi$ can be described as follows.
Let $C'\in\SVRT(w, \f_0)/\hspace{-.1cm}\sim$ and
 $T'\in C'$.
Let $T''$ be the Rothe
tableau defined by setting  $T''(B)=T'(B)$ if $B\in D(w)\setminus P(T',r+1)$, and setting
\begin{align*}
T''(B)=(T'(B)\setminus\{r,r+1\})\cup \{r\} \ \ \text{if $B\in P(T',r+1)$}.
\end{align*}
Notice that $T''\in\SVRT(w, \f_0)$.
We define $T$ as the Rothe tableau obtained from $T''$ by
adding the square $(r,w_r)$ filled with the set $\{r\}$, and then
 moving each square $B$ of $T''$ (together with
 the set $T''(B)$), which is to the right of the square $(r+1,w_r)$,
 up to row $r$.
By   construction, it is easily checked that $T\in\SVRT(ws_r, \f_0)$.
Let $C$ be the equivalence class in $\SVRT(ws_r, \f_0)/\hspace{-.1cm}\sim$
 containing $T$. Set $\Phi^{-1}(C')=C$.

Based on  Theorem \ref{mm-1} and Theorem \ref{mm-2},
we can establish the following relation.

\begin{theorem}\label{equivalence}
Let $w\neq w_0$ be a 1432-avoiding permutation, and $r$ be the first ascent of $w$.
For each equivalence class $C\in\SVRT(ws_r, \f_0)/\hspace{-.1cm}\sim$, we have
\begin{align}\label{F-3}
\pi_r G_{ws_r}(C; \x,\y)=G_w(\Phi(C); \x,\y).
\end{align}
\end{theorem}

The proof of  Theorem \ref{equivalence} requires  the
following property concerning  the operator $\pi_r$ due to Matsumura
\cite{Matsumura-2}.

\begin{lemma}[\mdseries{Matsumura \cite[Lemma 4.1]{Matsumura-2}}]\label{Matsu}
For an arbitrary  sequence $(a_1,a_2,\ldots, a_m)$ of positive integers,
\begin{align}
\pi_r((x_r\oplus y_{a_1})\cdots (x_r\oplus y_{a_m}))=&\sum_{k=1}^m
\ \prod_{j=1}^{k-1}(x_r\oplus y_{a_j})\prod_{j=k+1}^{m}(x_{r+1}\oplus y_{a_j})\nonumber\\[5pt]
& -\sum_{k=1}^{m-1}\
\prod_{j=1}^{k}(x_r\oplus y_{a_j})\prod_{j=k+1}^{m}(x_{r+1}\oplus y_{a_j}).
\label{ppqq-1}
\end{align}
Furthermore, the expression in  \eqref{ppqq-1} is symmetric in $x_r$ and $x_{r+1}$.
\end{lemma}

\noindent
 {\it Proof of Theorem \ref{equivalence}.}
 Assume that $T$ is any given Rothe tableau in $C$.
The polynomial  $h(C,i;\x,\y)$ defined in \eqref{hc}
has the following  reformulation:
\begin{align*}
h(C,i;\x,\y)=&\sum_{k=1}^{b_i(T)+1}\, \prod_{j=1}^{k-1} (-x_{r+1}\oplus y_{\ell_i(T)+j-1})
\prod_{j=k+1}^{b_i(T)+1} (-x_{r}\oplus y_{\ell_i(T)+j-1})\\[5pt]
&+\sum_{k=1}^{b_i(T)}\, \prod_{j=1}^k (-x_{r+1}\oplus y_{\ell_i(T)+j-1})
\prod_{j=k+1}^{b_i(T)+1} (-x_{r}\oplus y_{\ell_i(T)+j-1}).
\end{align*}
Hence $(-1)^{b_i(T)} h(C,i;\x,\y)$ coincides with
the right-hand side of \eqref{ppqq-1} by setting $m=b_i(T)+1$ and setting $a_j=\ell_i(T)+j-1$ for $1\leq j\leq m$, and then exchanging the variables
$x_r$ and $x_{r+1}$. It follows from Lemma \ref{Matsu}  that
\begin{align}\label{hpi}
h(C,i;\x,\y)&=(-1)^{b_i(T)}\cdot\pi_r
\left(\prod_{j=1}^{b_i(T)+1}(x_r\oplus y_{\ell_i(T)+j-1})\right)\nonumber\\
&=-\pi_r
\left(\prod_{j=1}^{b_i(T)+1}(-x_r\oplus y_{\ell_i(T)+j-1})\right),
\end{align}
which is a symmetric polynomial in
$x_r$ and $x_{r+1}$.

On the other hand, if a polynomial  $f(\x)$ is symmetric in $x_r$ and $x_{r+1}$, then for any polynomial $g(\x)$, it is easily checked that
\[\pi_r (f(\x) g(\x))= f(\x)\pi_r g(\x).\]
Therefore, applying $\pi_r$ to the formula of $G_{ws_r}(C; \x,\y)$ in Theorem \ref{mm-1},
we obtain that
\begin{align}\label{assd-1}
\pi_r G_{ws_r}(C; \x,\y)=&(-1)^{\ell(ws_r)} \left(\prod_{(i,j)\in D(ws_r)}\prod_{t\in T(i,j)\atop t\neq r,r+1}(-x_t\oplus y_{m_{ij}(ws_r)+i-t})\right)\nonumber\\[5pt]
&\  \cdot H_C(\x,\y)\cdot  J_C(\x,\y) \cdot \pi_r \left(\prod_{j=1}^{b_r(T)}(-x_r\oplus y_{\ell_r(T)+j})\right).
\end{align}

Let $T'\in\Phi(C)$ be any given Rothe tableau in the equivalent class of $\Phi(C)$.
By the construction of $\Phi$, it is easy to see  that
\begin{align}\label{assd-2}
\prod_{(i,j)\in D(w)}\prod_{t\in T(i,j)\atop t\neq r,r+1}(-x_t\oplus y_{m_{ij}(ws_r)+i-t})=\prod_{(i,j)\in D(w)}\prod_{t\in T'(i,j)\atop t\neq r,r+1}(-x_t\oplus y_{m_{ij}(w)+i-t}).
\end{align}
Again, by the construction of $\Phi$, it is also clear that for
$i>r+1$,
\[b_i(T)=b_i(T')\ \ \ \ \text{and}\ \ \ \ \ell_i(T)=\ell_i(T'),\]
which  imply that
\begin{align}
H_C(\x,\y)=H_{C'}(\x,\y).
\end{align}
Moreover, since $Q(T)=Q(T')$ and $m_{ij}(ws_r)=m_{ij}(w)$
for any $(i,j)\in Q^+(T)$, one has
\begin{align}
J_C(\x,\y)=J_{C'}(\x,\y).
\end{align}
Still, by the construction of $\Phi$, we see that
\[b_r(T)=b_{r+1}(T')+1\ \ \ \ \text{and}\ \ \ \ \ell_r(T)=\ell_{r+1}(T')-1.\]
So, by \eqref{hpi}, we have
\begin{align}\label{jfz}
\pi_r
\left(\prod_{j=1}^{b_{r}(T)}(-x_r\oplus y_{\ell_{r}(T)+j})\right)&=
\pi_r
\left(\prod_{j=1}^{b_{r+1}(T')+1}(-x_r\oplus y_{\ell_{r+1}(T')+j-1})\right)\nonumber\\[5pt]
&=-h(C',r+1;\x,\y).
\end{align}

Substituting  \eqref{assd-2}--\eqref{jfz} into \eqref{assd-1},
 we see that $\pi_r G_{ws_r}(C; \x,\y)=G_{w}(\Phi(C); \x,\y)$.
This completes the proof. \qed

By Theorem \ref{equivalence}  and  the bijection $\Phi$, we arrive at a proof of Theorem \ref{main-s}.  Using  induction on the length of $w$, we reach a proof of Theorem \ref{prop1}.

\subsection{Proof of Theorem \ref{prop2}}\label{sec4}

In this subsection, we confirm Theorem \ref{prop2}
by proving  the following statement.

\begin{theorem}\label{main3}
If  $w$ contains a 1432 pattern, then
\begin{align}\label{aognb}
\S_{w}(\x)\neq\sum_{T\in {\rm SRT}(w,\f_0)}
\prod_{(i,j)\in D(w)}\prod_{t\in T({i,j})}x_t.
\end{align}
\end{theorem}

By Corollary \ref{coro}, if $w$ is a 1432-avoiding permutation, then
$\mathfrak{S}_{w}(\x)$ must equal the right-hand side of \eqref{aognb}.
Hence Theorem \ref{main3} implies Theorem \ref{prop2}.

To finish the proof of Theorem \ref{main3},
we recall the balanced labeling model of Schubert polynomials
given by Fomin, Greene, Reiner
and Shimozono \cite{FoRe}.
To a square $(i,j)$ in the Rothe diagram $ D(w)$, the associated  hook  $H_{i,j}(w)$   is the collection of squares  $(i',j')$ of $D(w)$ such that either  $i'=i$ and $j'\geq j$, or $i'\geq i$ and $j'=j$.

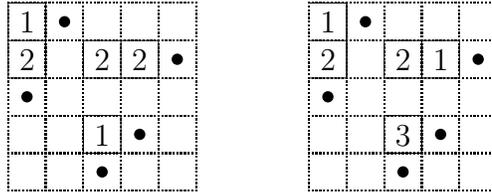
\begin{figure}[h]
\setlength{\unitlength}{0.5mm}
\begin{center}
\begin{picture}(300,55)

\qbezier[50](80,0)(105,0)(130,0)\qbezier[50](80,10)(105,10)(130,10)
\qbezier[50](80,20)(105,20)(130,20)\qbezier[50](80,30)(105,30)(130,30)
\qbezier[50](80,40)(105,40)(130,40)\qbezier[50](80,50)(105,50)(130,50)

\qbezier[50](80,0)(80,25)(80,50)\qbezier[50](90,0)(90,25)(90,50)
\qbezier[50](100,0)(100,25)(100,50)\qbezier[50](110,0)(110,25)(110,50)
\qbezier[50](120,0)(120,25)(120,50)\qbezier[50](130,0)(130,25)(130,50)

\put(95,45){\circle*{3}}\put(125,35){\circle*{3}}\put(85,25){\circle*{3}}
\put(115,15){\circle*{3}}\put(105,5){\circle*{3}}

\put(80,30){\line(1,0){10}}\put(80,40){\line(1,0){10}}\put(80,50){\line(1,0){10}}
\put(80,30){\line(0,1){20}}\put(90,30){\line(0,1){20}}

\put(100,30){\line(1,0){20}}\put(100,40){\line(1,0){20}}
\put(100,30){\line(0,1){10}}\put(110,30){\line(0,1){10}}\put(120,30){\line(0,1){10}}

\put(100,10){\line(1,0){10}}\put(100,20){\line(1,0){10}}
\put(100,10){\line(0,1){10}}\put(110,10){\line(0,1){10}}

\put(83,42){1}\put(83,32){2}\put(103,32){2}\put(113,32){2}
\put(103,12){1}


\qbezier[50](160,0)(185,0)(210,0)\qbezier[50](160,10)(185,10)(210,10)
\qbezier[50](160,20)(185,20)(210,20)\qbezier[50](160,30)(185,30)(210,30)
\qbezier[50](160,40)(185,40)(210,40)\qbezier[50](160,50)(185,50)(210,50)

\qbezier[50](160,0)(160,25)(160,50)\qbezier[50](170,0)(170,25)(170,50)
\qbezier[50](180,0)(180,25)(180,50)\qbezier[50](190,0)(190,25)(190,50)
\qbezier[50](200,0)(200,25)(200,50)\qbezier[50](210,0)(210,25)(210,50)

\put(175,45){\circle*{3}}\put(205,35){\circle*{3}}\put(165,25){\circle*{3}}
\put(195,15){\circle*{3}}\put(185,5){\circle*{3}}

\put(160,30){\line(1,0){10}}\put(160,40){\line(1,0){10}}\put(160,50){\line(1,0){10}}
\put(170,30){\line(0,1){20}}\put(170,30){\line(0,1){20}}

\put(180,30){\line(1,0){20}}\put(180,40){\line(1,0){20}}
\put(180,30){\line(0,1){10}}\put(190,30){\line(0,1){10}}\put(200,30){\line(0,1){10}}

\put(180,10){\line(1,0){10}}\put(180,20){\line(1,0){10}}
\put(180,10){\line(0,1){10}}\put(190,10){\line(0,1){10}}

\put(163,42){1}\put(163,32){2}\put(183,32){2}\put(193,32){1}
\put(183,12){3}

\end{picture}
\end{center}
\caption{Two  balanced  labelings for $w=25143$.}
\label{Rotheppp}
\end{figure}

A labeling $L$ of $D(w)$ is an assignment of positive
integers into the  squares of $D(w)$ such that each square
receives exactly one integer.   We use $L(i,j)$
to denote the label in the square $(i,j)\in D(w)$.
A labeling $L$  is called  balanced if for every square $(i,j)\in D(w)$, the label $L(i,j)$   remains unchanged after rearranging  the labels in the hook $H_{i,j}(w)$ so that they  are weakly increasing   from  right to left and from top to bottom.
Figure \ref{Rotheppp} illustrates two balanced labelings for the permutation $w=25143$.

A balanced  labeling of $D(w)$ is said to be column strict if
no column contains two equal labels. Let $\CSBL(w,\f_0)$  denote the set of
column strict balanced labelings of $D(w)$ such that $L(i,j)\leq i$ for each square $(i,j)\in D(w)$.
 Fomin, Greene, Reiner
and Shimozono \cite{FoRe} showed that
\begin{align}\label{Relate-Ba-2}
\mathfrak{S}_{w}(\x)=\sum_{L\in \CSBL(w,\f_0)} \prod_{(i,j)\in D(w)} x_{L(i,j)}.
\end{align}

We are now in a position to give a proof of Theorem \ref{main3}.

\noindent
{\it Proof of  Theorem \ref{main3}.} Assume that $w$  is a permutation
that contains a pattern 1432.
Recall that ${\rm SRT}(w,\f_0)$ is the set of single-valued Rothe tableaux of shape $D(w)$
flagged by $\f_0$.
By definition, it is clear that a Rothe tableau in $\SRT(w,\f_0)$ is a column strict balanced labeling, and hence  belongs to  $\CSBL(w,\f_0)$.
In view of \eqref{Relate-Ba-2}, to prove \eqref{aognb}, it suffices to show that there exits a balanced labeling in $\CSBL(w,\f_0)$ that does not belong to $\SRT(w,\f_0)$.
We next construct such a balanced  labeling $L$  in $\CSBL(w,\f_0)$.

Suppose that the subsequence $w_{i_1}w_{i_2}w_{i_3}w_{i_4}$ of $w$ has the same relative order as the pattern  1432, that is, $w_{i_1}<w_{i_4}<w_{i_3}<w_{i_2}$. Since $w_{i_3}>w_{i_4}$, there exists at least one square in the $i_3$-th row of $D(w)$. Let $(i_3,j)$ be the   rightmost  square in this row.
Let \[S=\{(i,j)\,|\,(i,j)\in D(w), i_1\leq i\leq i_3\}\]
be the subset of $D(w)$ consisting of the squares in column $j$ lying  between row $i_1$ and row $i_3$.
We  classify $S$ into two subsets according to whether a square $(i,j)\in S$ is the rightmost square in the row or not.
Let $S_1\subseteq S$ consists of square $(i,j)\in S$ such that $(i,j)$ is the rightmost square in row $i$. Clearly,  $S_1$ is nonempty since it  contains the square $(i_3,j)$. Let $S_2=S\backslash S_1$ be the complement. Since  $w_{i_2}>w_{i_3}$,
we see that the two squares $(i_2, j),(i_2, w_{i_3})$ belong to  $D(w)$. Hence $(i_2, j)\in S_2$, and so $S_2$ is also nonempty.

\begin{figure}[h]
\setlength{\unitlength}{0.5mm}
\begin{center}
\begin{picture}(250,105)
\qbezier[100](0,0)(50,0)(100,0)\qbezier[100](0,10)(50,10)(100,10)
\qbezier[100](0,20)(50,20)(100,20)\qbezier[100](0,30)(50,30)(100,30)
\qbezier[100](0,40)(50,40)(100,40)\qbezier[100](0,50)(50,50)(100,50)
\qbezier[100](0,60)(50,60)(100,60)\qbezier[100](0,70)(50,70)(100,70)
\qbezier[100](0,80)(50,80)(100,80)\qbezier[100](0,90)(50,90)(100,90)
\qbezier[100](0,100)(50,100)(100,100)

\qbezier[100](0,0)(0,50)(0,100)\qbezier[100](10,0)(10,50)(10,100)
\qbezier[100](20,0)(20,50)(20,100)\qbezier[100](30,0)(30,50)(30,100)
\qbezier[100](40,0)(40,50)(40,100)\qbezier[100](50,0)(50,50)(50,100)
\qbezier[100](60,0)(60,50)(60,100)\qbezier[100](70,0)(70,50)(70,100)
\qbezier[100](80,0)(80,50)(80,100)\qbezier[100](90,0)(90,50)(90,100)
\qbezier[100](100,0)(100,50)(100,100)

\put(5,95){\circle*{3}}\put(15,15){\circle*{3}}
\put(25,5){\circle*{3}}\put(35,85){\circle*{3}}
\put(45,75){\circle*{3}}\put(55,55){\circle*{3}}
\put(65,35){\circle*{3}}\put(75,25){\circle*{3}}
\put(85,65){\circle*{3}}\put(95,45){\circle*{3}}

\put(10,20){\line(1,0){20}}\put(10,30){\line(1,0){20}}
\put(10,40){\line(1,0){20}}\put(10,50){\line(1,0){20}}
\put(10,60){\line(1,0){20}}\put(10,70){\line(1,0){20}}
\put(10,80){\line(1,0){20}}\put(10,90){\line(1,0){20}}

\put(10,20){\line(0,1){70}}\put(20,20){\line(0,1){70}}
\put(30,20){\line(0,1){70}}

\put(50,60){\line(1,0){30}}\put(50,70){\line(1,0){30}}

\put(50,60){\line(0,1){10}}\put(60,60){\line(0,1){10}}
\put(70,60){\line(0,1){10}}\put(80,60){\line(0,1){10}}

\put(60,40){\line(1,0){20}}\put(60,50){\line(1,0){20}}

\put(60,40){\line(0,1){10}}\put(70,40){\line(0,1){10}}
\put(80,40){\line(0,1){10}}

\put(-7,93){\text{\small $i_1$}}\put(-7,43){\text{\small $i_2$}}
\put(-7,33){\text{\small $i_3$}}\put(-7,13){\text{\small $i_4$}}

\put(22, 83){{\small $\spadesuit$}}\put(22, 73){{\small $\spadesuit$}}
\put(22, 63){{\small $\clubsuit$}}\put(22, 53){{\small $\spadesuit$}}
\put(22, 43){{\small $\clubsuit$}}\put(22, 33){{\small $\spadesuit$}}

\put(45,-10){\small (a)}


\qbezier[100](150,0)(200,0)(250,0)\qbezier[100](150,10)(200,10)(250,10)
\qbezier[100](150,20)(200,20)(250,20)\qbezier[100](150,30)(200,30)(250,30)
\qbezier[100](150,40)(200,40)(250,40)\qbezier[100](150,50)(200,50)(250,50)
\qbezier[100](150,60)(200,60)(250,60)\qbezier[100](150,70)(200,70)(250,70)
\qbezier[100](150,80)(200,80)(250,80)\qbezier[100](150,90)(200,90)(250,90)
\qbezier[100](150,100)(200,100)(250,100)

\qbezier[100](150,0)(150,50)(150,100)\qbezier[100](160,0)(160,50)(160,100)
\qbezier[100](170,0)(170,50)(170,100)\qbezier[100](180,0)(180,50)(180,100)
\qbezier[100](190,0)(190,50)(190,100)\qbezier[100](200,0)(200,50)(200,100)
\qbezier[100](210,0)(210,50)(210,100)\qbezier[100](220,0)(220,50)(220,100)
\qbezier[100](230,0)(230,50)(230,100)\qbezier[100](240,0)(240,50)(240,100)
\qbezier[100](250,0)(250,50)(250,100)

\put(155,95){\circle*{3}}\put(165,15){\circle*{3}}
\put(175,5){\circle*{3}}\put(185,85){\circle*{3}}
\put(195,75){\circle*{3}}\put(205,55){\circle*{3}}
\put(215,35){\circle*{3}}\put(225,25){\circle*{3}}
\put(235,65){\circle*{3}}\put(245,45){\circle*{3}}

\put(160,20){\line(1,0){20}}\put(160,30){\line(1,0){20}}
\put(160,40){\line(1,0){20}}\put(160,50){\line(1,0){20}}
\put(160,60){\line(1,0){20}}\put(160,70){\line(1,0){20}}
\put(160,80){\line(1,0){20}}\put(160,90){\line(1,0){20}}

\put(160,20){\line(0,1){70}}\put(170,20){\line(0,1){70}}
\put(180,20){\line(0,1){70}}

\put(200,60){\line(1,0){30}}\put(200,70){\line(1,0){30}}

\put(200,60){\line(0,1){10}}\put(210,60){\line(0,1){10}}
\put(220,60){\line(0,1){10}}\put(230,60){\line(0,1){10}}

\put(210,40){\line(1,0){20}}\put(210,50){\line(1,0){20}}

\put(210,40){\line(0,1){10}}\put(220,40){\line(0,1){10}}
\put(230,40){\line(0,1){10}}

\put(143,93){\text{\small $i_1$}}\put(143,43){\text{\small $i_2$}}
\put(143,33){\text{\small $i_3$}}\put(143,13){\text{\small $i_4$}}
\put(143,53){\text{\small $i_0$}}

\put(172, 83){{\small $\diamondsuit$}}\put(172, 73){{\small $\diamondsuit$}}
\put(172, 53){{\small $\diamondsuit$}}

\put(195,-10){\small (b)}

\end{picture}
\end{center}
\caption{An example for the proof of Theorem \ref{main3}.}
\label{Exam}
\end{figure}
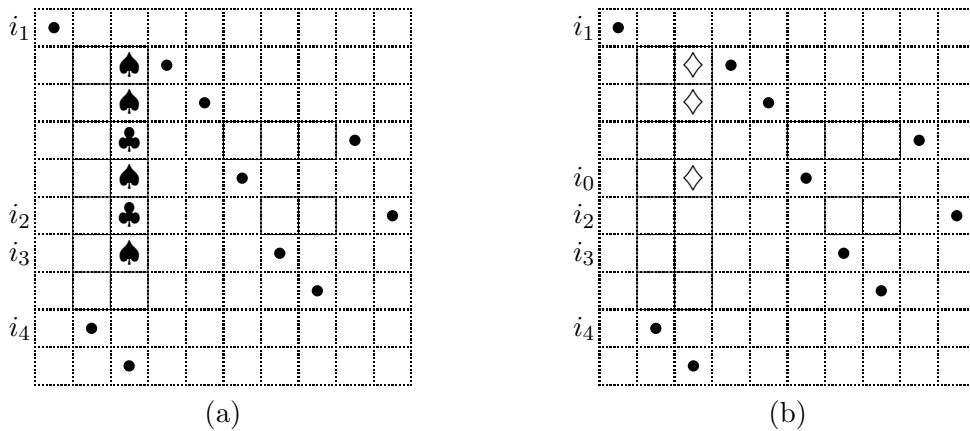

Let us use an example in Figure \ref{Exam} to illustrate the sets $S_1$ and $S_2$.
In this example, $w=1\,4\,5\,9\,6\,10\,7\,8\,2\,3$ and  the  subsequence $w_1w_6w_7w_9$  forms a 1432-pattern.
The rightmost square of $D(w)$ in the $i_3$-th row is the square $(7,3)$, and so we have
\[S=\{(i,3)\,|\, i=2,3,4,5,6,7\}.\]
Moreover, the squares belonging to $S_1$ and $S_2$ are signified with
$\spadesuit$ and $\clubsuit$ in Figure \ref{Exam}(a), respectively.

Let $i_0$ be the  smallest row index such that:
(1) the square  $(i_0,j)\in S_1$; (2) there exists a square in $S_2$ lying  above $(i_0,j)$.
Such an row index exists since the
$i_3$-th row satisfies the above  conditions.
Let $S'=\{(i,j)\in S_1\,|\,i\le i_0\}$ be the subset of $S_1$ including the squares above $(i_0,j)$.
In the example in Figure \ref{Exam}, we see that $i_0=5$ and the squares of $S'$ are signified with the symbol  $\diamondsuit$.

Assume that $|S'|=k$ and $(r_1,j), \ldots, (r_k,j)$ are the squares of $S'$, where
$r_1<\cdots<r_k=i_0$. Note  that $r_1>i_1$. This is because  $w_{i_1}$ is the smallest element of $\{w_{i_1},w_{i_2},w_{i_3},w_{i_4}\}$ and thus the square $(i_1,j)\notin D(w)$.

We now construct a balanced   labeling $L$ of $D(w)$ as follows. If a square $(s,t)$ of $D(w)$ is not contained in $S'$, then we set $L(s,t)=s$. For the squares $(r_1,j), \ldots, (r_k,j)$  of $S'$,  we set $L(r_1,j)=i_1$  and $L(r_{p},j)=r_{p-1}$ for $p=2,\ldots,k$.
For the permutation in Figure
\ref{Exam}, the labeling $L$ is given in Figure \ref{Exams}, where the integers in $S'$ are written in boldface.

\begin{figure}[h]
\setlength{\unitlength}{0.5mm}
\begin{center}
\begin{picture}(100,105)
\qbezier[100](0,0)(50,0)(100,0)\qbezier[100](0,10)(50,10)(100,10)
\qbezier[100](0,20)(50,20)(100,20)\qbezier[100](0,30)(50,30)(100,30)
\qbezier[100](0,40)(50,40)(100,40)\qbezier[100](0,50)(50,50)(100,50)
\qbezier[100](0,60)(50,60)(100,60)\qbezier[100](0,70)(50,70)(100,70)
\qbezier[100](0,80)(50,80)(100,80)\qbezier[100](0,90)(50,90)(100,90)
\qbezier[100](0,100)(50,100)(100,100)

\qbezier[100](0,0)(0,50)(0,100)\qbezier[100](10,0)(10,50)(10,100)
\qbezier[100](20,0)(20,50)(20,100)\qbezier[100](30,0)(30,50)(30,100)
\qbezier[100](40,0)(40,50)(40,100)\qbezier[100](50,0)(50,50)(50,100)
\qbezier[100](60,0)(60,50)(60,100)\qbezier[100](70,0)(70,50)(70,100)
\qbezier[100](80,0)(80,50)(80,100)\qbezier[100](90,0)(90,50)(90,100)
\qbezier[100](100,0)(100,50)(100,100)

\put(5,95){\circle*{3}}\put(15,15){\circle*{3}}
\put(25,5){\circle*{3}}\put(35,85){\circle*{3}}
\put(45,75){\circle*{3}}\put(55,55){\circle*{3}}
\put(65,35){\circle*{3}}\put(75,25){\circle*{3}}
\put(85,65){\circle*{3}}\put(95,45){\circle*{3}}

\put(10,20){\line(1,0){20}}\put(10,30){\line(1,0){20}}
\put(10,40){\line(1,0){20}}\put(10,50){\line(1,0){20}}
\put(10,60){\line(1,0){20}}\put(10,70){\line(1,0){20}}
\put(10,80){\line(1,0){20}}\put(10,90){\line(1,0){20}}

\put(10,20){\line(0,1){70}}\put(20,20){\line(0,1){70}}
\put(30,20){\line(0,1){70}}

\put(50,60){\line(1,0){30}}\put(50,70){\line(1,0){30}}

\put(50,60){\line(0,1){10}}\put(60,60){\line(0,1){10}}
\put(70,60){\line(0,1){10}}\put(80,60){\line(0,1){10}}

\put(60,40){\line(1,0){20}}\put(60,50){\line(1,0){20}}

\put(60,40){\line(0,1){10}}\put(70,40){\line(0,1){10}}
\put(80,40){\line(0,1){10}}

\put(-7,93){\text{\small $i_1$}}\put(-7,43){\text{\small $i_2$}}
\put(-7,33){\text{\small $i_3$}}\put(-7,13){\text{\small $i_4$}}

\put(13,82){2}\put(13,72){3}
\put(13,62){4}\put(13,52){5}
\put(13,42){6}\put(13,32){7}
\put(13,22){8}

\put(23,82){{\bf 1}}\put(23,72){{\bf 2}}
\put(23,62){4}\put(23,52){{\bf 3}}
\put(23,42){6}\put(23,32){7}
\put(23,22){8}

\put(53,62){4}\put(63,62){4}\put(73,62){4}

\put(63,42){6}\put(73,42){6}


\end{picture}
\end{center}
\caption{A balanced labeling in $\CSBL(w, \f_0)$, but not in $\SRT(w, \f_0)$.}
\label{Exams}
\end{figure}

By the construction of $L$, it is not hard to check that
$L$ is a column strict balanced labeling in $\CSBL(w, \f_0)$.
Moreover,  the entries in the $j$-th column of $L$ are not increasing.
So $L$ does not belong to $\SRT(w,\f_0)$. This completes the proof.
\qed

\section{Rothe tableau complexes}\label{sec5}

In this section, we prove the tableau formulas  in Theorem \ref{3-1}.
To do this, we investigate  the structure of  Rothe tableau complexes, which  is a specific family of the tableau complexes
as introduced by Knutson, Miller and Yong \cite{KnMiYo-2}.
Using Theorem \ref{main} and the properties of  tableau complexes established in \cite{KnMiYo-2}, we
obtain  two alternative tableau formulas for the Grothendieck polynomials of 1432-avoiding permutations, as given in Theorem \ref{3-1}.

Let us proceed with a brief  review of the   Hilbert series  of the
Stanley-Reisner ring of a simplicial complex,  see   \cite{MiSt,Stan} for more detailed information.
An (abstract) simplicial complex $\Delta$ on a finite   vertex set $V$
is a collection of subsets of $V$ such that if $\sigma\in \Delta$ and $\tau\subseteq \sigma$, then $\tau\in \Delta$.
Each subset  $\sigma\in \Delta$ is  called a face of $\Delta$.
A face $\sigma$ is called a facet of $\Delta$ if $\sigma$ is not a subset of any other faces.
Clearly,  $\Delta$ is  determined by its   facets.

Let $\mathds{k}[\mathbf{t}]$ be the ring of polynomials over a field
$\mathds{k}$ in the variables $t_v$ where $v\in V$. The Stanley-Reisner ideal $I_{\Delta}$
is the ideal generated by the monomials corresponding to the subsets  of $V$
that are not faces of $\Delta$,
namely,
\[I_\Delta=\left<\prod_{v\in \tau}t_v\,|\, \text{$\tau\subseteq V$, but $\tau\not\in \Delta$}\right>. \]
The Stanley-Reisner ring  of $\Delta$, denoted $\mathds{k}[\Delta]$,
is the quotient ring $\mathds{k}[\mathbf{t}]/ I_\Delta$.
The Hilbert series $H(\mathds{k}[\Delta]; \mathbf{t})$ of $\mathds{k}[\mathbf{t}]/ I_\Delta$ is equal to the sum of
monomials not belonging to $I_{\Delta}$.
It is well known \cite{MiSt,Stan} that  $H(\mathds{k}[\Delta]; \mathbf{t})$
has the following formula:
\[H(\mathds{k}[\Delta]; \mathbf{t})=
\frac{K(\mathds{k}[\Delta];\mathbf{t})}{\prod_{v\in V}(1-t_v)},\]
where
\[K(\mathds{k}[\Delta];\mathbf{t})=\sum_{\sigma\in \Delta}\prod_{v\in \sigma}t_v\prod_{v\not\in \sigma }(1-t_v).\]
The numerator $K(\mathds{k}[\Delta];\mathbf{t})$ is called the $K$-polynomial
of $\mathds{k}[\Delta]$.

We now restrict attention to  the $K$-polynomials  of
 tableau complexes introduced in  \cite{KnMiYo-2}.
Let $X$ and $Y$ be two finite sets.
A map $f$ from $X$ to $Y$ is called a tableau, which can be viewed
 as an assignment of  elements of $Y$ to elements of $X$ such that each $x\in X$
receives exactly one element of $Y$.
A tableau $f$ can also be identified with the following set
\[\{(x\mapsto y)\,|\, \text{$x\in X$ and $f(x)=y$}\}\subseteq X\times Y\]
 of ordered pairs.  Let $U$ be a subset of tableaux
from $X$ to $Y$, and let $E\subseteq X\times Y$ be a set of ordered pairs  such that $f\subseteq E$ for each $f\in U$. The tableau
complex corresponding to $U$ and $E$, denoted  $\Delta_E(X\xrightarrow{U}Y)$, can be defined as follows. Let us first define
a simplex $\Delta_E$. For each pair $(x\mapsto a)\in E$, write  $(x\arrownot\mapsto y)=E\setminus \{(x\mapsto y)\}$ for the complement of $\{(x\mapsto y)\}$, and let
\[V=\{(x\arrownot\mapsto y)\,|\, (x\mapsto y)\in E\}.\]
Denote by $\Delta_E$  the simplex with vertex set $V$, that is,
$\Delta_E$ is the collection of all of the subsets of $V$.

Let $F \subseteq V$ be a face of  $\Delta_E$. Assume that $F$ has $k$ vertices $  (x_1\arrownot\mapsto y_1),\ldots, (x_k\arrownot\mapsto y_k)$.
Then $F$ can be identified with the following subset of $E$:
\[E\setminus \{(x_i \mapsto y_i)\,|\, 1\leq i\leq k\}.\]
On the other hand, each subset of $E$ can be viewed as a
 set-valued tableau  from
$X$ to $Y$, that is, a map that assigns each element of $X$
with a subset   of $Y$.  To be more specific,
for a subset $A$ of $E$, the corresponding set-valued tableau
is defined by assigning $x\in X$ with the subset $\{y\in Y\,|\,
(x  \mapsto y)\in A\}$. So the face $F$ of $\Delta_E$
can also be identified with a set-valued tableau such that
for $x\in X$,
\[F(x)=\{y\in Y\,|\, (x  \mapsto y)\in E,  \text{$(x  \mapsto y)\neq (x_i \mapsto y_i)$ for $1\leq i\leq k$}\}.\]
From now on, a face $F$ of $\Delta_E$ can  be identified either with
a subset of $E$ or with a set-valued tableau from $X$ to $Y$,
which will not cause confusion  from the context.
By the definition of $\Delta_E$,  a vertex $(x\arrownot\mapsto y)\in V$
belongs to $F$ if  and only if the pair $(x\mapsto y) $ does not belong to $F$.

Recall that $U$ is a set of tableaux from $X$ to $Y$
such that $f\subseteq E$ for each  $f\in U$. So each tableau $f$
in $U$
is a face of $\Delta_E$.
The tableau complex $\Delta_E(X\xrightarrow{U}Y)$ is defined as the
subcomplex of $\Delta_E$
such that the  facets of $\Delta_E(X\xrightarrow{U}Y)$ are the tableaux   in $U$.
This means that a set-valued tableau $F\subseteq E$
is a face of $\Delta_E(X\xrightarrow{U}Y)$ if and only if $F$ contains some
tableau $f\in U$.

When $X$ and $Y$ are further endowed with  partially ordered structures,
Knutson, Miller and Yong \cite{KnMiYo-2} found three different  expressions for
the $K$-polynomial of a tableau complex.

\begin{theorem}[\mdseries{Knutson-Miller-Yong \cite{KnMiYo-2}}]\label{KMY}
Let $X$ and $Y$ be two finite posets. For each $x\in X$, let $Y_x$ be a totally
ordered subset of $Y$.
Let $\Psi$ be a set of pairs $(x,x')$ in
$X$ with $x<x'$. Let $U$ be the set of tableaux $f\colon X\rightarrow Y$
such that
\begin{itemize}
\item[(a)] $f(x)\in Y_x$;

\item[(b)] $f$ is weakly order preserving, that is, if $x\leq x'$,
then $f(x)\leq f(x')$;

\item[(c)] if $(x,x')\in \Psi$, then $f(x)<f(x')$.
\end{itemize}
Set $E=\bigcup_{f\in U}f$. Let $\mathbf{t}=\{t_{(x\arrownot\mapsto a)}\,|\, (x\arrownot\mapsto a)\in V\}$ be the set of variables corresponding to
the vertices of the tableau complex $\Delta=\Delta_E(X\xrightarrow{U}Y)$. Then,
$\Delta$ is homeomorphic to a ball or a sphere. Moreover, the corresponding $K$-polynomial
has the following expressions.
\begin{itemize}

    \item[1.] Let $U_1$ be the set of set-valued tableaux $F\subseteq E$ such that  every tableau $f\subseteq F$ lies in $U$. Then,
    \begin{equation}\label{KMY-1}
    K(\mathds{k}[\Delta];\mathbf{t})=\sum_{F\in U_1} (-1)^{|F|-|X|} \prod_{x\in X}\prod_{a\in F(x)}\left(1-t_{(x\arrownot\mapsto a)}\right).
    \end{equation}

\item[2.] Let $U_2$ be the set of  set-valued tableaux $F\subseteq E$ each
    containing some tableau $f\in U$. Then,
\begin{equation}\label{KMY-2}
K(\mathds{k}[\Delta];\mathbf{t})=\sum_{F\in U_2} \prod_{x\in X}\left(\prod_{a\in F(x)}\left(1-t_{(x\arrownot\mapsto a)}\right)\prod_{a\in E(x)\setminus F(x)}t_{(x\arrownot\mapsto a)}
\right).
    \end{equation}

    \item[3.] Given a tableau $f\in U$ and $x\in X$, let $Y_f(x)$ be the set
    of $y\in Y$ such that $f(x)< y$ and moving the label on $x$ from $f(x)$
    up to $y$ still yields a tableau in $U$. Then,
    \begin{equation}\label{KMY-3}
    K(\mathds{k}[\Delta];\mathbf{t})=\sum_{f\in U}
     \prod_{x\in X}\left(\left(1-t_{(x\arrownot\mapsto f(x))}\right)\prod_{a\in Y_f(x)}t_{x\arrownot\mapsto a}\right).
     \end{equation}
\end{itemize}

\end{theorem}

We now  consider the specific tableau  complex such that the facets are the
  single-valued Rothe tableaux in  $\SRT(w, \f_0)$. To be consistent with the aforementioned  notation, let $X=D(w)$ and   $Y$  be the set of positive
 integers. Set $U=\SRT(w, \f_0)$
and
\[E=\bigcup _{T\in \SRT(w, \f_0)} T.\]
We denote the above defined tableau complex by $\Delta(w)=\Delta_E(X\xrightarrow{U}Y)$,
and call $\Delta(w)$ the Rothe tableau complex for $w$.

Using Theorem \ref{main} and Theorem \ref{KMY}, we can now give a proof of Theorem \ref{3-1}.

\noindent{\it Proof of Theorem \ref{3-1}}. We  define a partial ordering on $D(w)$ as follows.
For two distinct squares $B$ and $B'$ of $D(w)$, we use $B\rightarrow B'$ to represent
that either $B$ and $B'$  are in the same row and $B$ lies to the right of $B'$,
or $B$ and $B'$  are in the same column and $B$ lies above $B'$.
Define  $B<B'$ if there exists a sequence
 $(B=B_1,B_2,\ldots,B_k=B')$ of squares of $D(w)$ such that
 \[B=B_1\rightarrow B_2\rightarrow \cdots \rightarrow B_k=B'.\]

For each square $B=(i,j)$ of $D(w)$, let $Y_B=\{1,2,\ldots,i\}$.
Moreover, we set $\Psi$ to be the set of pairs  $(B,B')$ with $B<B'$
such that $B$ and $B'$ are in the same column of $D(w)$.
Now we see that the tableaux satisfying the conditions $(a)$,
$(b)$ and $(c)$ in Theorem \ref{KMY} are exactly the single-valued Rothe tableaux in $\SRT(w,\f_0)$.  Recall that the set $U_1$
defined in Theorem \ref{KMY}  consists of the set-valued tableaux $F\subseteq E$ such that  every tableau in $F$ lies in $U$. Clearly,  $F\subseteq E$ is a set-valued tableau  satisfying  that  every tableau contained in $F$ lies in $U$ if and only if
$F$ is a set-valued Rothe tableau in  $\SVRT(w,\f_0)$.
Thus we have $U_1=\SVRT(w, \f_0)$. Replacing  $t_{x\arrownot\mapsto a}$ with $x=(i,j)\in D(w )$ by
\[\frac{x_a}{y_{m_{ij}(w)+i-a}}\]
 and then replacing $x_t$ by $1-x_t$ and
$y_t$ by $\frac{1}{1-y_t}$, the $K$-polynomial $  K(\mathds{k}[\Delta];\mathbf{t})$
in \eqref{KMY-1} becomes
\begin{align*}
  \sum_{T\in \SVRT(w, \f_0)} (-1)^{|T|-\ell(w)} \prod_{(i,j)\in D(w)}\prod_{t\in T(i,j)}
  (x_t\oplus y_{m_{ij}(w)+i-t}),
    \end{align*}
which agrees with the formula \eqref{CFG-R}
in Theorem \ref{main}. Making the same substitutions
in \eqref{KMY-2} and \eqref{KMY-3}, we are led to
\eqref{3-1-1} and \eqref{3-1-2}
respectively. This completes the proof.
\qed

\vspace{0.5cm}
 \noindent{\bf Acknowledgments.} This work was supported by  the 973
Project and the National Science
Foundation of China.

\end{document}